\documentclass[a4paper,12pt]{article}
\pdfoutput=1
\usepackage{etex}

\usepackage{amssymb}
\usepackage{amsmath}
\usepackage{amsthm}
\usepackage{bbm}

\usepackage{lineno}
\usepackage[cm]{fullpage}
\usepackage{graphicx}
\usepackage{wrapfig}
\usepackage{subfig}
\usepackage{float}
\usepackage{booktabs}
\usepackage{caption}
\usepackage[inline]{asymptote}

\allowdisplaybreaks

\makeatletter
\newsavebox{\sfe@box}
\newenvironment{subfloatenv}[1]{%
\def\sfe@caption{#1}%
\setbox\sfe@box\hbox\bgroup\color@setgroup}%
{\color@endgroup\egroup\subfloat[\sfe@caption]%
{\usebox{\sfe@box}}}
\makeatother

\newcommand{\norm}[1]{\ensuremath{\lVert#1\rVert}}

\newcommand{\reverse}[1]{\ensuremath{\widetilde{#1}}}

\newcommand{\R}[1]{\ensuremath{\mathbb{R}^{#1}}}

\newcommand{\E}[1]{\ensuremath{\mathbb{E}^{#1}}}
\newcommand{\M}[1]{\ensuremath{\mathbb{M}^{#1}}}

\newcommand{\tb}[1]{\ensuremath{\textbf{#1}}}
\newcommand{\mb}[1]{\ensuremath{\boldsymbol{#1}}}

\newcommand{\e}{\tb{e}} % basis vectors
\newcommand{\I}{\tb{I}} % pseudoscalar
\title{Clifford algebra and the projective model of Minkowski (pseudo-Euclidean) spaces}
\author{Andrey Sokolov}

\begin{document}
\maketitle

\begin{abstract}
I apply the algebraic framework developed in  \cite{gunn2011geometry} to Minkowski (pseudo-Euclidean) spaces in 2, 3, and 4 dimensions.
The exposition follows the template established in arXiv:1307.2917 \cite{sokolov2013clifford} for Euclidean spaces.
The emphasis is on geometric structures, but some contact with special relativity is made by considering relativistic addition of velocities
and Lorentz transformations, both of which can be seen as rotation applied to points and to lines.
The language used in the paper reflects the emphasis on geometry, rather than applications to special relativity.
The use of Clifford algebra greatly simplifies the study of Minkowski spaces, since unintuitive synthetic techniques are
replaced by algebraic calculations.
\end{abstract}

\tableofcontents

\section{Minkowski plane \M{2} ( $1+1$ spacetime)}

To simplify the comparison with the 4-dimensional Minkowski space, I will replace \(y\) with \(t\) and \(b\) with \(h\), so that
\(\tb{a}=d\e_0+a\e_1+h\e_2\) and \(\tb{P}=w\e_{12}+x\e_{20}+t\e_{01}\).
In special relativity, \(t\) is interpreted as a time coordinate and \(x\) as a spatial coordinate.
A finite line \(\tb{a}\) in \M{2} is called proper if \(\tb{a}^2>0\), improper if \(\tb{a}^2<0\), and null if  \(\tb{a}^2=0\).
As in Euclidean space, \(\e_0^2=0\) for the line at infinity in \M{2}, but \(\e_0\) is not called null as it is not finite.
For a line \(\tb{a}\),  I get \(\tb{a}^2=a^2-h^2\) and its norm is given by \(\norm{\tb{a}}=|a^2-h^2|^{\tfrac{1}{2}}\).
The norm of a finite point is the same as in Euclidean space, i.e.\ \(\norm{\tb{P}}=|w|\).
Note, however, that \(\tb{P}^2=w^2\) in \M{2}, whereas \(\tb{P}^2=-w^2\) in \E{2}.

In Euclidean space, the terms \emph{finite} and \emph{at infinity} match the metric properties of the corresponding lines and points.
This is no longer the case in non-Euclidean spaces and in Minkowski spaces in particular.
Nevertheless, due to the fact that Minkowski space is flat, it resembles Euclidean space more than other non-Euclidean spaces
and it is still convenient to continue to use the same terminology even though it no longer precisely matches the metric properties.
The term \emph{at infinity} applied to points in Minkowski space is a shorthand for describing points represented by stacks of lines rather than sheaves.
Note that if a point in Euclidean space is at an infinite distance from the origin, 
then it is either at an infinite distance in Minkowski space or the distance to it is undefined.
The line at infinity consists of points represented by stacks only (\(d\e_0\) is the corresponding vector in \R{3*}); any other line will be called finite.

\begin{figure}[t!]
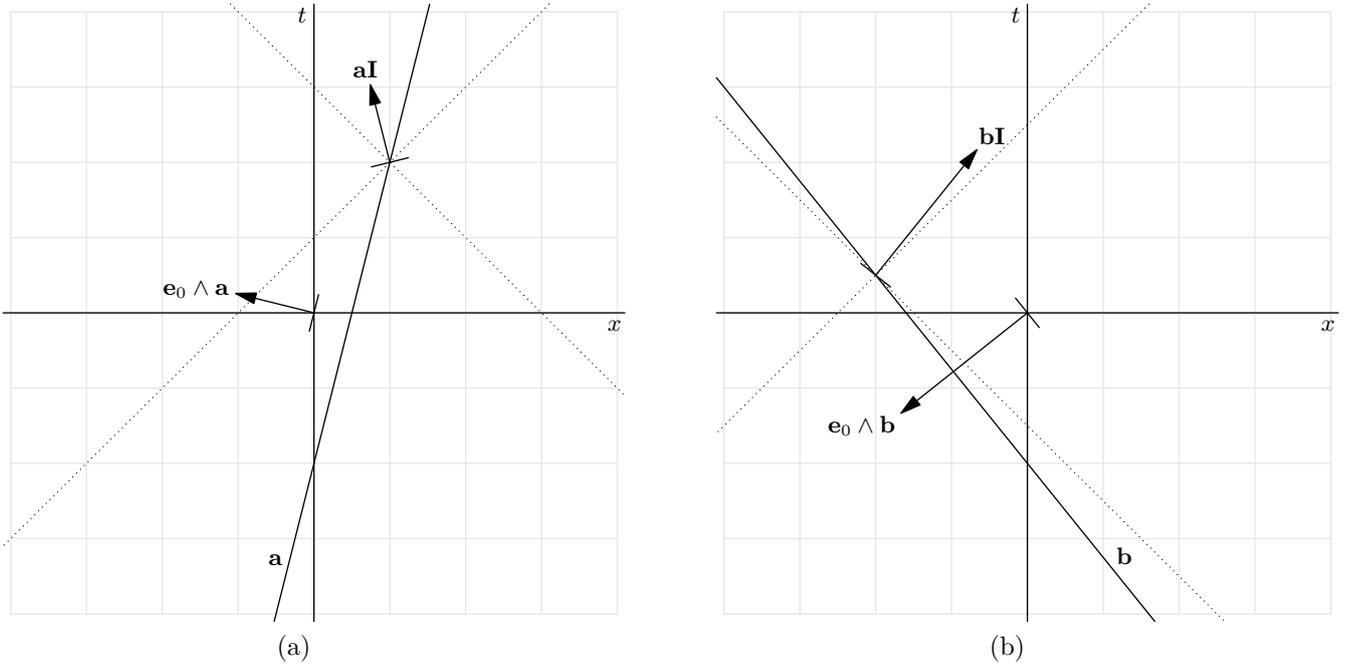

\begin{subfloatenv}{ }
\begin{asy}
import Figure2D;
Figure f = Figure(yaxis_name="$t$");
metric = Metric(Minkowski);

pair centre = (1,2);
f.null_lines(centre);

MV a = -0.25*join(Point(1,0,-2),Point(1,1,2));
a=a/norm(a);
f.line(a,label="$\textbf{a}$",position=0.1,align=W,draw_orientation=false);

MV N = wedge(e_0,a);
f.point_at_infinity(N, "$\textbf{e}_0\wedge\textbf{a}$");
f.point_at_infinity(a*I, "$\textbf{aI}$", O=centre);

\end{asy}
\end{subfloatenv}\hfill%
\begin{subfloatenv}{  }
\begin{asy}
import Figure2D;
Figure f = Figure(yaxis_name="$t$");
metric = Metric(Minkowski);
pair centre = (-2,0.5);
f.null_lines(centre);

MV a = -0.5*join(Point(1,0,-2),Point(1,-2,0.5));
a=a/norm(a);
f.line(a,label="$\textbf{b}$",position=0.1,draw_orientation=false);

MV N = wedge(e_0,a);
f.point_at_infinity(N, "$\textbf{e}_0\wedge\textbf{b}$");
f.point_at_infinity(a*I, "$\textbf{bI}$", O=centre);

\end{asy}
\end{subfloatenv}
\caption{Orientation and polar points in \M{2}}
\label{orientation and polar points in M2}

\end{figure}

Any structure or relation in \M{2} that does not depend on the metric, such as intersection of lines, is the same as in Euclidean space \E{2}.
For example, the point at infinity \(\e_0\wedge\tb{a}\) lies on \(\tb{a}\) and provides the top-down orientation of the line.
Recall that in the depiction of points at infinity such as \(\e_0\wedge\tb{a}\) (see Figure~\ref{orientation and polar points in M2}), 
the arrow specifies the top-down orientation of the point at infinity.
In Euclidean space, it plays an additional role since the direction of the arrow is perpendicular to the line.
The notion of perpendicularity depends on the metric, so the above observation does not necessarily apply in \M{2}.
A point at infinity which lies in the direction perpendicular to \(\tb{a}\) is given by the polar point \(\tb{a}\I\).
This is illustrated in Figure~\ref{orientation and polar points in M2} for two proper lines \(\tb{a}\) and \(\tb{b}\),
 where relevant null lines are shown with dotted lines.

\begin{figure}[h!]
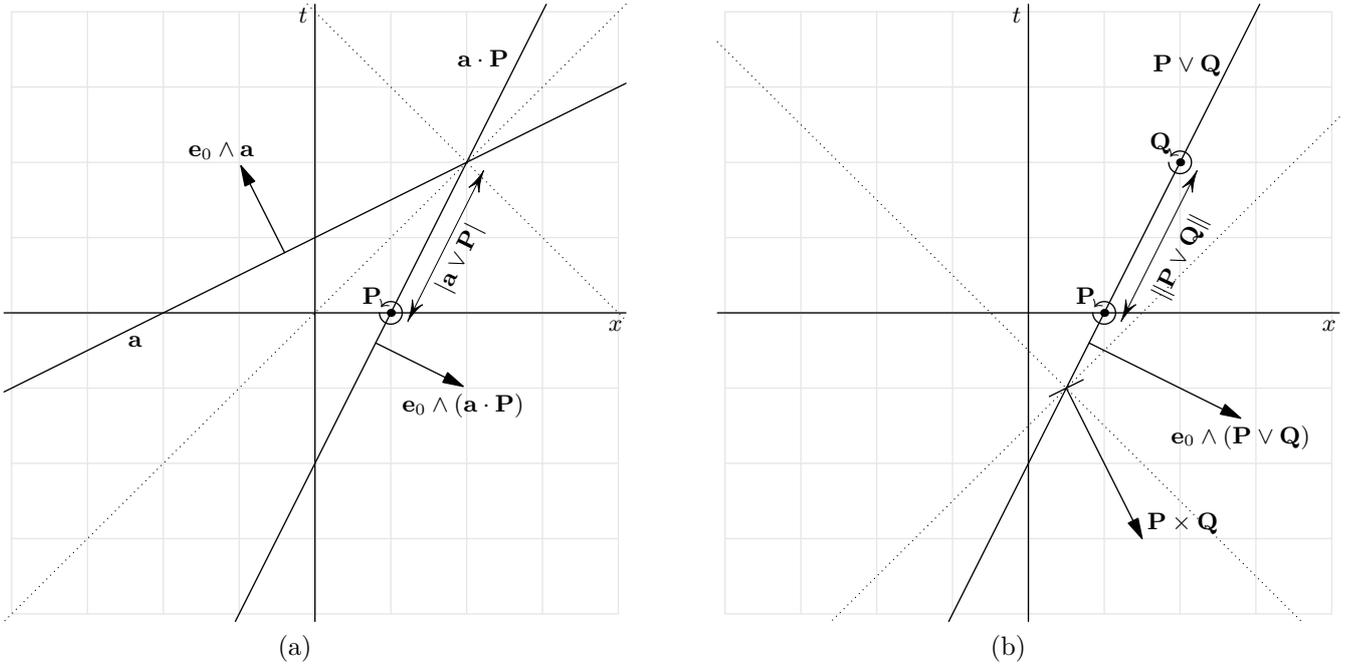

\begin{subfloatenv}{ }
\begin{asy}
import Figure2D;
metric = Metric(Minkowski);
Figure f = Figure(yaxis_name="$t$");

var a = -join(Point(1,0,1),Point(1,2,2));
a=a/norm(a);
var P = Point(1,1,0);

f.line(a, "$\textbf{a}$", "$\textbf{e}_0\wedge\textbf{a}$");
f.point(P, "$\textbf{P}$",align=NW);

var b = dot(a,P);
f.line(b, "$\textbf{a}\cdot\textbf{P}$", position=0.1, "$\textbf{e}_0\wedge(\textbf{a}\cdot\textbf{P})$", orientation_align=S);
f.stretch(P, dot(P,a)/a, "$|\textbf{a}\vee\textbf{P}|$",label_shift=(-10,0),align=E,switch_shift_direction=true);

f.null_lines(topair(wedge(a,b)));

\end{asy}
\end{subfloatenv}\hfill%
\begin{subfloatenv}{ }
\begin{asy}
import Figure2D;
metric = Metric(Minkowski);
Figure f = Figure(yaxis_name="$t$");

pair centre = (0.5,-1);
f.null_lines(centre);

MV P = Point(1, 1, 0);
MV Q = Point(1, 2, 2);

f.point(P, "\textbf{P}",align=NW);
f.point(Q, "\textbf{Q}",align=NW);
f.line(join(P,Q), "$\textbf{P}\vee\textbf{Q}$", position=0.1, align=W, "$\textbf{e}_0\wedge(\textbf{P}\vee\textbf{Q})$",orientation_align=S);
f.point_at_infinity(cross(P,Q),"$\textbf{P}\times\textbf{Q}$", align=NE, O=centre);
f.stretch(P, Q, "$\norm{\textbf{P}\vee\textbf{Q}}$",align=E,switch_shift_direction=true);

\end{asy}
\end{subfloatenv}
\caption{Points and lines in \M{2}}
\label{points and lines in M2}
\end{figure}

Perpendicularity is further illustrated in Figure~\ref{points and lines in M2}(a), where the line \(\tb{a}\cdot\tb{P}\)
is perpendicular to \(\tb{a}\).
On the other hand, parallel lines are those lines that intersect at infinity, so their mutual appearance in \M{2} is the same as in \E{2}.

The join of two points \(\tb{P}\) and  \(\tb{Q}\) (see Figure~\ref{points and lines in M2}(b))
is independent of the metric and, therefore, the same as in \E{2}.
The commutator \(\tb{P}\times\tb{Q}\) is a point at infinity in the direction perpendicular to the line \(\tb{P}\vee\tb{Q}\).
The distance \(r\) between normalised points  \(\tb{P}\) and \(\tb{Q}\) connected by a proper line, i.e.\ \((\tb{P}\vee\tb{Q})^2>0\),
is given by   
\begin{equation}
r=\norm{\tb{P}\vee\tb{Q}}.
\label{distance in M2}
\end{equation}
This expression is equivalent to the one used in Euclidean space because both Euclidean and Minkowski spaces have a degenerate distance measure.
The norm depends on the metric, so the value of distance in Minkowski space is different from that in Euclidean space.
Substituting  \(\tb{P}=\e_{12}+x_P\e_{20}+t_P\e_{01}\)  and \(\tb{Q}=\e_{12}+x_Q\e_{20}+t_Q\e_{01}\) yields 
the familiar formula \(r=\sqrt{(t_P-t_Q)^2-(x_P-x_Q)^2}\).
If \(\tb{P}\vee\tb{Q}\) is null or improper, the distance is undefined.
Note that if the expression~(\ref{distance in M2}) is applied to points connected by a null line, then one formally obtains \(r=0\)
which may be useful in some applications, but I do not extend the definition of distance to such cases.
It is possible to formally apply (\ref{distance in M2}) to points connected by null or improper lines, but 
the restriction of the definition of distance to points connected by proper lines only is necessary
in order to enjoy the properties expected of distance such as additivity.
Indeed, if a point \(\tb{S}\) lies between \(\tb{P}\) and \(\tb{Q}\) on the line \(\tb{P}\vee\tb{Q}\), which is proper, 
then \(r=r_{SP}+r_{SQ}\), where \(r_{SP}\) and \(r_{SQ}\) are the distances from \(\tb{S}\) to \(\tb{P}\) and \(\tb{Q}\) respectively.
In special relativity, the points \(\tb{P}\) and  \(\tb{Q}\) are treated as events and the distance between then
is usually called the spacetime interval.
The distance between a normalised line \(\tb{a}\) and a normalised point \(\tb{P}\) is given by
\(|\tb{a}\vee\tb{P}|\) if the line \(\tb{a}\) is improper, which implies \(\tb{a}\cdot\tb{P}\) is proper (see Figure~\ref{points and lines in M2}(a)).

\begin{figure}[t!]
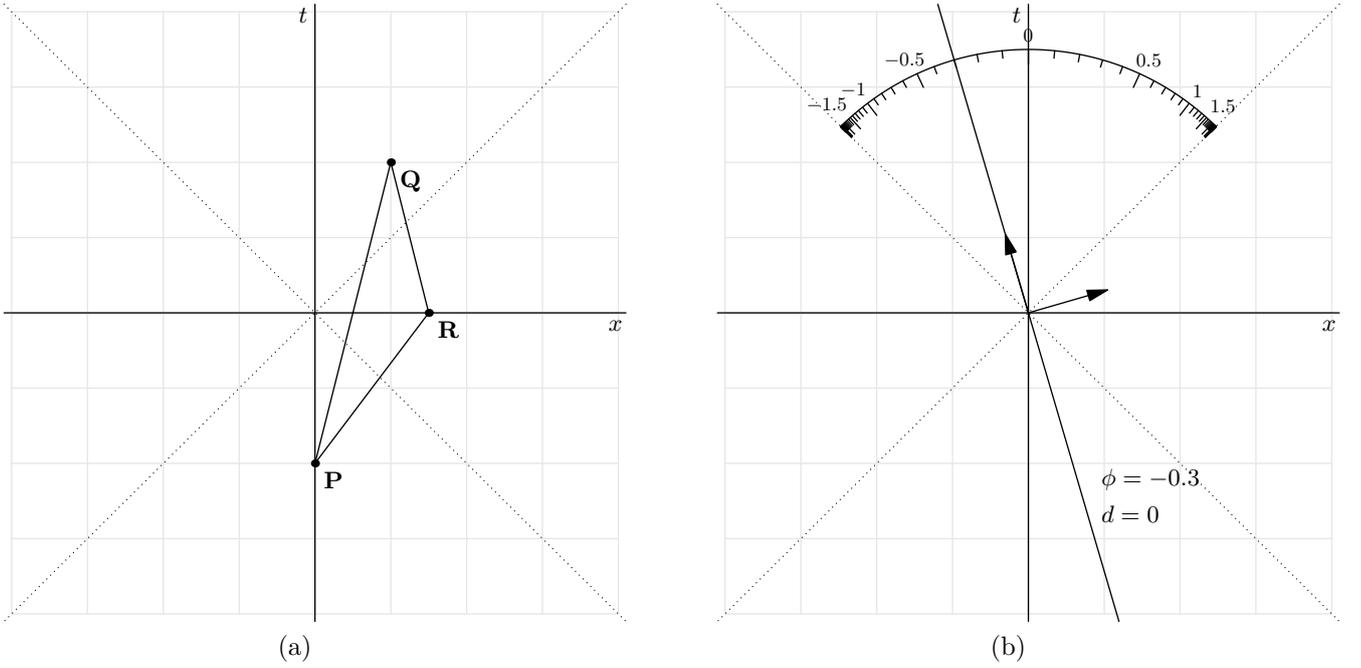

\begin{subfloatenv}{ }
\begin{asy}
import Figure2D;
Figure f = Figure(yaxis_name="$t$");
metric = Metric(Minkowski);
f.null_lines();

Point P = Point(1,0,-2);
Point Q = Point(1,1,2);
Point R = Point(1,1.5,0);

f.point(P,"$\textbf{P}$",draw_orientation=false);
f.point(Q,"$\textbf{Q}$",draw_orientation=false);
f.point(R,"$\textbf{R}$",draw_orientation=false);

f.interval(P,Q);
f.interval(Q,R);
f.interval(R,P);

\end{asy}
\end{subfloatenv}\hfill%
\begin{subfloatenv}{  }
\begin{asy}
import Figure2D;
Figure f = Figure(yaxis_name="$t$");
metric = Metric(Minkowski);
f.null_lines();

path p = Arc((0,0),3.5,45,45+90);
real f(real phi) { pair v = (sinh(phi),cosh(phi)); v=5*v/length(v); path h = (0,0)--v; return intersect(p,h)[0]; }

axis(p, LeftTicks(scale(0.8)*Label("%3.1g"), OmitFormat(-2,-2.5,-3,-3.5,-4,-4.5,-5,2,2.5,3,3.5,4,4.5,5), N=20,n=5), ticklocate(-5,5,f), above=true);

real b = -0.3;
MV a = Line(0,cosh(b),-sinh(b)); 
f.line(a, label="$\begin{aligned}&\phi=-0.3\\ &d=0\end{aligned}$",position=0.8,align=(2,0), draw_bottom_up_orientation=true);

\end{asy}
\end{subfloatenv}
\caption{Angles and distances in \M{2}}
\label{angles and distances in M2}
\end{figure}

For normalised points \(\tb{P}\),  \(\tb{Q}\) and \(\tb{R}\)  connected by proper lines as in Figure~\ref{angles and distances in M2}(a), 
the following inequality is satisfied
\begin{equation}
\norm{\tb{P}\vee\tb{Q}}\ge \norm{\tb{P}\vee\tb{R}} + \norm{\tb{R}\vee\tb{Q}},
\end{equation}
where the sign of inequality is the opposite of what it would be in Euclidean space.
So, a straight line provides the longest distance between two points; the distance of any other proper path is shorter.
This is particularly obvious when  \(\tb{P}\) and \(\tb{Q}\) are connected by a vertical line, so that
\(\norm{\tb{P}\vee\tb{Q}}=|t_P-t_Q|\), and  \(\tb{R}\) is such that the lines \(\tb{P}\vee\tb{R}\) and \(\tb{R}\vee\tb{Q}\)
are very close to null lines, so that \(\norm{\tb{P}\vee\tb{R}}\ll1\) and \(\norm{\tb{R}\vee\tb{Q}}\ll1\).
The area \(\cal S\) of a triangle defined by points \(\tb{P}\), \(\tb{Q}\) and \(\tb{R}\)
is given by 
\begin{equation}
{\cal S}=\frac{1}{2!}|\tb{P}\vee\tb{Q}\vee\tb{R}|
\end{equation}
if the points are normalised and connected by proper lines.
The area is the same as in \E{2}, since \(\tb{P}\vee\tb{Q}\vee\tb{R}\) is a scalar, the normalisation of points is the same as in Euclidean space,
and the join does not depend on the metric.
%The area is undefined if the lines connecting \(\tb{P}\), \(\tb{Q}\) and \(\tb{R}\) are not all proper.

In Euclidean space, any finite line can be parametrised by the angle the line makes with a certain fixed direction
(usually the \(x\)-axis in the direction of increasing \(x\)).
In Minkowski space, proper lines can be parametrised in a similar fashion. 
For instance, any normalised proper line whose bottom-up orientation is facing in the direction of increasing \(t\), 
i.e.\ towards the future in the context of special relativity,
 can be written as
\begin{equation}
\tb{a}=d\e_0 + \e_1\cosh{\phi} - \e_2\sinh{\phi},
\label{line parameter phi}
\end{equation}
where \(\phi\in(-\infty,+\infty)\) is the angle the line makes with the vertical direction as illustrated in Figure~\ref{angles and distances in M2}(b).
I will refer to such lines as future-oriented for brevity.
A proper line \(\tb{a}\) is future-oriented if \(\tb{a}\cdot\e_1>0\), where \(\e_1\) is the \(t\)-axis (\(\e_1\) is future-oriented by this definition).
The top-down orientation of the future-oriented lines is facing in the direction of increasing \(x\) (see Figure~\ref{angles and distances in M2}(b)).
The parametrisation (\ref{line parameter phi}) is not applicable to null or improper lines. 
Nor is it applicable to proper lines whose bottom-up orientation is facing in the direction of decreasing \(t\), which will be called past-oriented.
A proper line \(\tb{a}\) is past-oriented if \(\tb{a}\cdot\e_1<0\). 
If \(\tb{a}\) is past-oriented, then \(-\tb{a}\) is future-oriented, so the past-oriented lines can be parametrised by
\begin{equation}
\tb{a}=-(d\e_0 + \e_1\cosh{\phi} - \e_2\sinh{\phi}),
\label{line parameter phi 2}
\end{equation}
where \(\phi\) is the same as in (\ref{line parameter phi}).

In the following, no distinction is made between future- and past-oriented lines, unless stated otherwise.
No physical significance is assigned to this classification.
The separation of future- and past-oriented lines is geometric in nature.
For instance, any future-oriented line can be obtained from \(\e_1\) by a proper motion and
 any past-oriented line can be obtained from \(-\e_1\).
However, there are no proper motions that can convert a future-oriented line into a past-oriented line, and vice versa.

A proper line represents the  worldline of an object moving at a constant velocity \(u=\tanh{\phi}\),
assuming the units where the speed of light \(c=1\).
The angle \(\phi\) can be thought of as a parametrisation of relativistic velocities.
It is sometimes referred to as rapidity, but it has a simple geometric meaning as the angle between the worldline it corresponds to 
and the \(t\)-axis.

\begin{figure}[t!]
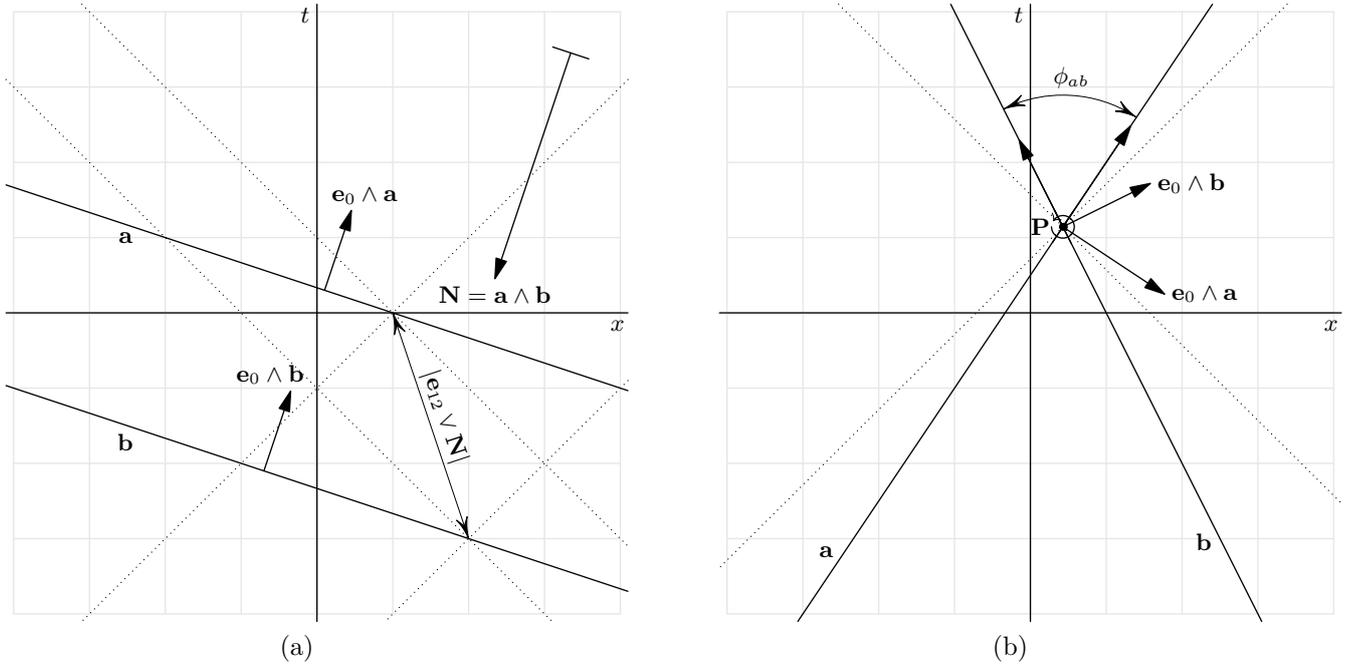

\begin{subfloatenv}{ }
\begin{asy}
import Figure2D;
Figure f = Figure(yaxis_name="$t$");
metric = Metric(Minkowski);

MV a = Line(-1,1,3);
a = a/norm(a);
MV Na = wedge(e_0,a);

MV b = Line(7,1,3);
b = b/norm(b);
MV Nb = wedge(e_0,b);

MV N = wedge(a,b);

f.point_at_infinity(N, O=(3.34655,3.4535), "$\textbf{N}=\textbf{a}\wedge\textbf{b}$", align=S);
f.line(a, "$\textbf{a}$", "$\textbf{e}_0\wedge\textbf{a}$");
f.line(b, "$\textbf{b}$", "$\textbf{e}_0\wedge\textbf{b}$", orientation_align=(-0.5,1));

MV a_ = join(Point(1,2,-3), a*I);
f.stretch(wedge(b,a_), wedge(a,a_), "$|\textbf{e}_{12}\vee\textbf{N}|$", shift_distance=0.0,label_shift=20.0,align=(0.5,1.5),rotation_angle=180);
f.null_lines(wedge(b,a_));
f.null_lines(wedge(a,a_));

\end{asy}
\end{subfloatenv}\hfill%
\begin{subfloatenv}{ }
\begin{asy}
import Figure2D;
Figure f = Figure(yaxis_name="$t$");
metric = Metric(Minkowski);

MV a =  Line(1,3,-2);
a = a/norm(a);
MV Na = wedge(e_0,a);

MV b = -Line(2,-2,-1);
b = b/norm(b);
MV Nb = wedge(e_0,b);

MV P = wedge(a,b);
pair O = topair(P);

f.point(P, "$\textbf{P}$", align=W);
f.line(a, "$\textbf{a}$", position=0.9,draw_orientation=false,draw_bottom_up_orientation=true, O=O);
f.point_at_infinity(Na, O=O, "$\textbf{e}_0\wedge\textbf{a}$",align=E);
f.line(b, "$\textbf{b}$", position=0.86,draw_orientation=false,draw_bottom_up_orientation=true, O=O);
f.point_at_infinity(Nb, O=O, "$\textbf{e}_0\wedge\textbf{b}$", align=E);
f.arc(P, -Na,-Nb,"$\phi_{ab}$", N, rotation=90,radius=1.75);
f.null_lines(P);

\end{asy}
\end{subfloatenv}
\caption{Parallel and intersecting lines in \M{2}}
\label{parallel and intersecting lines}
\end{figure}

As in Euclidean space, \(\tb{a}\tb{b}=\pm1+\tb{a}\wedge\tb{b}\) for two normalised parallel (anti-parallel) 
lines \(\tb{a}\) and \(\tb{b}\),
where \(\tb{a}\wedge\tb{b}\) is a point at infinity where the lines intersect. % (\(-1\) for anti-parallel lines). 
The distance between \(\tb{a}\) and \(\tb{b}\)  can be computed 
with \(\norm{\e_{12}\vee(\tb{a}\wedge\tb{b})}\),
provided that the lines are normalised  and improper (see Figure~\ref{parallel and intersecting lines}(a)).
If \(\tb{a}\) and \(\tb{b}\) are proper or null, the distance between them is undefined.

For two normalised proper lines \(\tb{a}\) and \(\tb{b}\) representing worldlines as discussed above, the inner product satisfies \(|\tb{a}\cdot\tb{b}|\ge1\).
Furthermore, if both lines are future-oriented or past-oriented, then \(\tb{a}\cdot\tb{b}\ge1\).
In this case, 
the geometric product is given by
\begin{equation}
\tb{a}\tb{b}=\cosh{\phi_{ab}}+\tb{P}\sinh{\phi_{ab}},
\label{angle in M2}
\end{equation}
where \(\tb{P}\) is a normalised point at the  intersection of the lines and
\(\phi_{ab}\in(0,+\infty)\) is the angle between the lines (see Figure~\ref{parallel and intersecting lines}(b))
defined by
\begin{equation}
\cosh\phi_{ab}=\tb{a}\cdot\tb{b}.
\end{equation}
%The point \(\tb{P}\) is counterclockwise by definition, which determines the sign of \(\phi_{ab}\), i.e.\ 
%\(\phi_{ab}>0\) if \(\tb{a}\wedge\tb{b}\) is counterclockwise and \(\phi_{ab}<0\) if \(\tb{a}\wedge\tb{b}\) is clockwise
%(\(\phi_{ab}=0\) if the lines coincide).
%The expression \(\cosh{\phi_{ab}} = \tb{a}\cdot\tb{b}\) does not alone define the angle since both \(\phi_{ab}\) and \(-\phi_{ab}\) satisfy it.
The angle is undefined if \(\tb{a}\) is future-oriented and \(\tb{b}\) is past-oriented, or vice versa.
The existence of the angle between two lines implies that there is a proper motion that converts one line into the other.
Since this is impossible if \(\tb{a}\) is future-oriented and \(\tb{b}\) is past-oriented, or vice versa, 
I do not define the angle between such lines.
Like distance, the angle defined by (\ref{angle in M2}) is additive in \M{2}.

Recall that in Euclidean space the angle between the lines satisfies \(\cos\alpha=\tb{a}\cdot\tb{b}\),
whereas in \M{2} hyperbolic functions are used in the definition.
This discrepancy is due to the fact that the angular measure in Euclidean space is elliptic and in Minkowski space it is hyperbolic.
The angular measure is hyperbolic in the other kinematic spaces, de-Sitter and anti de-Sitter, as well.

If \(\tb{a}\cdot\tb{b}\ge1\),
then the angle between \(\tb{a}\) and \(\tb{b}\) is defined and \(\phi_{ab}=|\phi_a-\phi_b|\), where
\(\tb{a}\) is parametrised with the angle \(\phi_a\) and \(\tb{b}\) with \(\phi_b\).
This is illustrated in Figure~\ref{parallel and intersecting lines}(b) where
\(\tb{a}= \tfrac{1}{\sqrt{5}}(\e_0+3\e_1-2\e_2)\) and
\(\tb{b}= \tfrac{1}{\sqrt{3}}(-2\e_0+2\e_1+\e_2)\) are normalised.
The parametrisation angles satisfy \(\tanh{\phi_a}=\tfrac{2}{3}\) and 
 \(\tanh{\phi_b}=-\tfrac{1}{2}\) 
, which gives \(\phi_a\approx0.80\), \(\phi_b\approx-0.55\), and \(\phi_{ab}\approx1.35\).
Using the addition formula for hyperbolic tangent gives \(\tanh{\phi_{ab}}=\tfrac{7}{8}\).

\begin{figure}[t!]
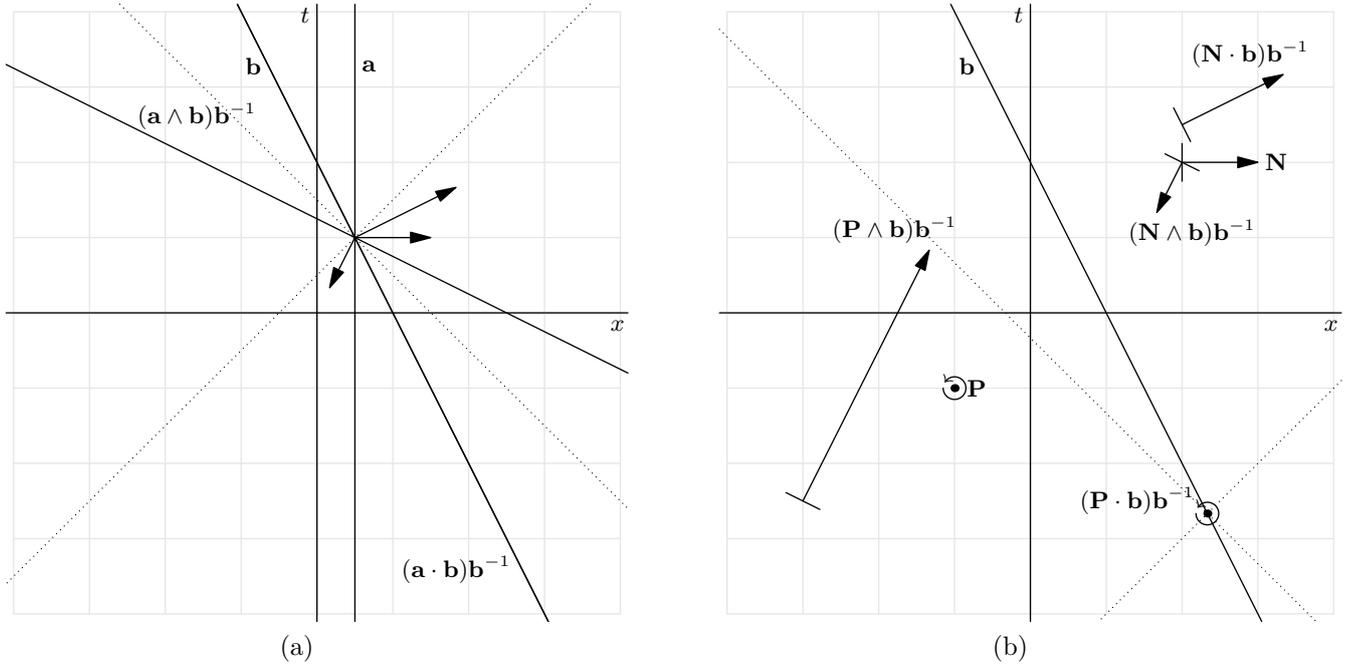

\begin{subfloatenv}{ }
\begin{asy}
import Figure2D;
Figure f = Figure(yaxis_name="$t$");
metric = Metric(Minkowski);

var a = -Line(1/2,-1,0);
var b = -Line(2,-2,-1)/4;

f.line(b, "$\textbf{b}$", position=0.1,align=W, draw_orientation=false);

f.line(a,"$\textbf{a}$", position=0.1, align=E,draw_orientation=true, O=(1/2,1));
f.line(dot(a,b)/b, "$(\textbf{a}\cdot\textbf{b})\textbf{b}^{-1}$", position=0.9, draw_orientation=true, O=(1/2,1));
f.line(wedge(a,b)/b, "$(\textbf{a}\wedge\textbf{b})\textbf{b}^{-1}$", position=0.8, align=(1,0.5),draw_orientation=true, O=(1/2,1));

f.null_lines(wedge(a,b));

\end{asy}
\end{subfloatenv}\hfill%
\begin{subfloatenv}{ }
\begin{asy}
import Figure2D;
Figure f = Figure(yaxis_name="$t$");
metric = Metric(Minkowski);

var b = Line(2,-2,-1)/4;
var P = Point(1,-1,-1);

f.line(b, "$\textbf{b}$", position=0.9,align=W, draw_orientation=false);
f.point(P, "$\textbf{P}$", align=E);

var Q = dot(P,b)/b;
f.point(Q, "$(\textbf{P}\cdot\textbf{b})\textbf{b}^{-1}$", align=(-1,0.5));
f.point_at_infinity(wedge(P,b)/b,"$(\textbf{P}\wedge\textbf{b})\textbf{b}^{-1}$",O=(-3,-2.5),align=(-0.5,1));

var N = wedge(e_0,Line(0,1,0));
f.point_at_infinity(N, "$\textbf{N}$", O=(2,2));
f.point_at_infinity(dot(N,b)/b, "$(\textbf{N}\cdot\textbf{b})\textbf{b}^{-1}$", O=(2,2.5),align=(-0.5,1));
f.point_at_infinity(wedge(N,b)/b, "$(\textbf{N}\wedge\textbf{b})\textbf{b}^{-1}$", O=(2,2), align=(0.5,-1));

f.null_lines(dot(P,b)/b);

\end{asy}
\end{subfloatenv}
\caption{Projection and rejection in \M{2}}
\label{projection and rejection in M2}
\end{figure}

Projection, rejection, reflection, and scaling in \M{2} are defined in the same way as in \E{2}.
Projection on and rejection by a line in \M{2} is illustrated in Figure~\ref{projection and rejection in M2}.
Projection on and rejection by a point in \M{2} is identical to that in \E{2}.

\begin{figure}[t!]
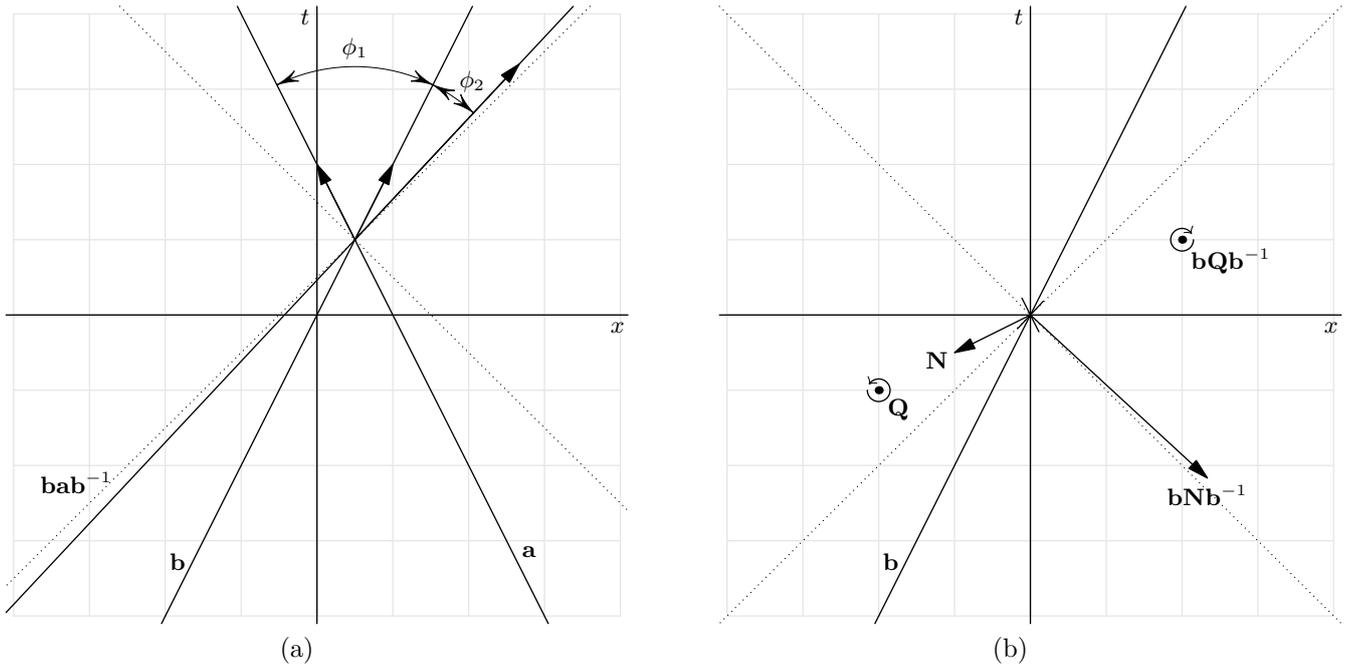

\begin{subfloatenv}{ }
\begin{asy}
import Figure2D;
Figure f = Figure(yaxis_name="$t$");
metric = Metric(Minkowski);

var a = -Line(1,-1,-1/2);
var b = 2*Line(0,1/2,-1/4);
var P = wedge(a,b);

write(b);
f.line(a, "$\textbf{a}$", position=0.9, align=NE, draw_orientation=false, O=topair(P), draw_bottom_up_orientation=true);
f.line(b, "$\textbf{b}$", position=0.9, align=W, draw_orientation=false, O=topair(P), draw_bottom_up_orientation=true);

f.line(b*a/b, "$\textbf{b}\textbf{a}\textbf{b}^{-1}$", draw_orientation=false,O=topair(P), position=0.8, align=(-1,0.5), draw_bottom_up_orientation=true);

f.arc(P, -wedge(e_0,b), -wedge(e_0,a),"$\phi_1$", N, rotation=90,radius=2.3);
f.arc(P, -wedge(e_0,b*a/b), -wedge(e_0,b),"$\phi_2$", NE, rotation=90,radius=2.3);
f.null_lines(P);

\end{asy}
\end{subfloatenv}\hfill%
\begin{subfloatenv}{ }
\begin{asy}
import Figure2D;
Figure f = Figure(yaxis_name="$t$");
metric = Metric(Minkowski);

var a = Line(1,-1,-1/2);
var b = Line(0,1/2,-1/4);
var P = wedge(a,b);

f.line(b, "$\textbf{b}$", position=0.9, align=W, draw_orientation=false);

var N = wedge(e_0, Line(1,-1,-1/2));
f.point_at_infinity(N, "$\textbf{N}$", O=(0,0));
f.point_at_infinity(b*N/b, "$\textbf{b}\textbf{N}\textbf{b}^{-1}$", align=(0,-1), O=(0,0));

var Q = Point(1,-2,-1);
f.point(Q,"$\textbf{Q}$");
f.point(b*Q/b,"$\textbf{b}\textbf{Q}\textbf{b}^{-1}$");
f.null_lines();
\end{asy}
\end{subfloatenv}
\caption{Reflection in a line in \M{2}}
\label{reflections in a line in M2}
\end{figure}

I adopt the same expressions for the top-down and bottom-up reflections as those used in Euclidean space.
A bottom-up reflection in a line is illustrated in Figure~\ref{reflections in a line in M2}.
For the bottom-up reflection in a proper line \(\tb{b}\), if \(\tb{a}\) is future-oriented, then 
\(\tb{b}\tb{a}\tb{b}^{-1}\) is also future-oriented and \(\phi_1=\phi_2\),
where \(\phi_1\) is the angle between \(\tb{b}\) and \(\tb{a}\) and \(\phi_2\) is the angle between \(\tb{b}\) and \(\tb{b}\tb{a}\tb{b}^{-1}\)
(\(\tb{b}\) is also future-oriented, so the angles are defined).
Moreover, the angle between \(\tb{a}\) and \(\tb{b}\tb{a}\tb{b}^{-1}\) is given by \(\phi_1+\phi_2\) due to additivity.
Reflection in a point is the same as in \E{2}.

The Spin group in \M{2} is defined in the same way as in \E{2}.
Namely, it consists of multivectors that can be written as 
the product of an even number of normalised proper lines, i.e.\ lines that square to unity.
Spinors are even and each spinor satisfies \(S\reverse{S}=1\).
Any proper motion in \M{2} can be generated by a spinor that has the form
\(S=e^A\), where \(A\) is a bivector.
Spinors in \M{2} form a Lie group, whose Lie algebra consists of bivectors.
In Minkowski space, the Lie group consists of multivectors \(e^{A}\) and \(-e^A\), where \(A\) is an arbitrary bivector.
For example, the spinor \(S=\tb{a}\tb{b}\tb{c}\tb{d}\), where
\(\tb{a}=\tfrac{1}{4}\e_0+\sqrt{2}\e_1+\e_2\), 
\(\tb{b}=\sqrt{2}\e_1+\e_2\), 
\(\tb{c}=\tfrac{1}{4}\e_0-\sqrt{2}\e_1+\e_2\), and
\(\tb{d}=-\sqrt{2}\e_1+\e_2\) are all proper and normalised, 
 can be written as \(S=e^{-\tfrac{1}{2}\e_{20}}\).
There are spinors in \M{2} that can only be written as \(-e^A\), e.g.\ for \(S=\tb{a}\tb{b}\) where \(\tb{a}=-\e_1\) and \(\tb{b}=\e_1\),
I have \(-e^0=-\e_1\e_1=-1\) and \(e^A\ne-1\) for any bivector \(A\) (incidentally, \(e^{\pi\e_{12}}=-1\) in \E{2}).
Note that the action of \(S=-e^A\) is identical to that of \(e^{A}\), since the minus sign  cancels out in \(SMS^{-1}\).

\begin{figure}[t!]
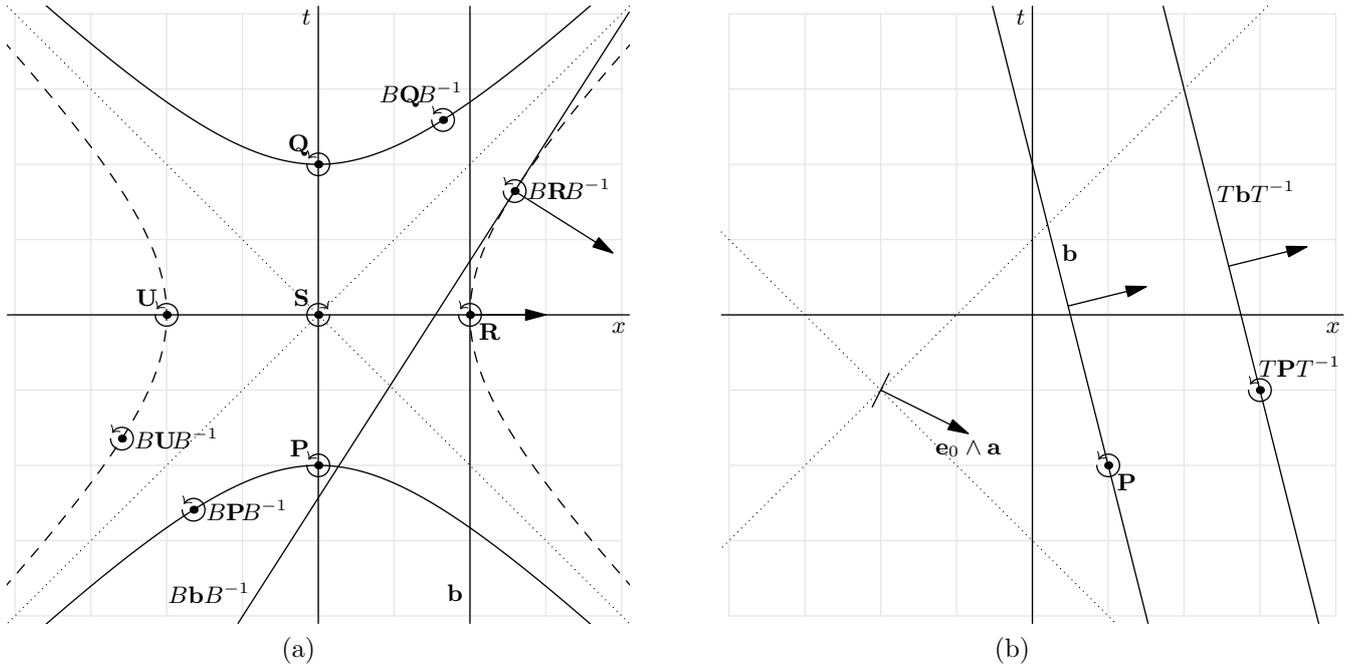

\begin{subfloatenv}{ }
\begin{asy}
import Figure2D;
Figure f = Figure(yaxis_name="$t$");
metric = Metric(Minkowski);

var S = -Point(1,0,0);
f.point(S, "$\textbf{S}$", draw_orientation=true,align=NW);

f.null_lines(S);

MV b = -Line(1,-1/2,0);
b=b/norm(b);

real phi1 = -2; real phi2 = 2; pair c = topair(S); real x0 = c.x; real t0 = c.y; real r = 2;
pen p = defaultpen; // Dotted

real x(real phi) {return x0+r*sinh(phi);}
real t(real phi) {return  t0+r*cosh(phi);}
draw(graph(x,t,phi1,phi2),p);

real x(real phi) {return x0-r*sinh(phi);}
real t(real phi) {return  t0-r*cosh(phi);}
draw(graph(x,t,phi1,phi2),p);

f.crop();

bool show_tangent_circle = true;
if (show_tangent_circle) { 
  real rt = 2;
  real enl = 1.5;
  real x(real phi) {return x0+rt*cosh(phi);}
  real t(real phi) {return  t0+rt*sinh(phi);}
  draw(graph(x,t,enl*phi1,enl*phi2),p=dashed);

  real x(real phi) {return x0-rt*cosh(phi);}
  real t(real phi) {return  t0-rt*sinh(phi);}   
  draw(graph(x,t,enl*phi1,enl*phi2),p=dashed);
  f.crop();
}

var Q = Point(1,0,2);
var P = Point(1,0,-2);
var R = Point(1,2,0);
var U = Point(1,-2,0);
f.point(Q, "$\textbf{Q}$",align=NW);
f.point(P, "$\textbf{P}$",align=NW);
f.point(U, "$\textbf{U}$",align=NW);
f.point(R, "$\textbf{R}$");

real phi = 0.75;
var Rot = exp(-1/2*phi*S);

f.point(Rot*R/Rot, "$B\textbf{R}\!B^{-1}$",align=(1,0.1));
f.point(Rot*Q/Rot, "$B\textbf{Q}\!B^{-1}$",align=(-0.4,1.));
f.point(Rot*P/Rot, "$B\textbf{P}\!B^{-1}$",align=E);
f.point(Rot*U/Rot, "$B\textbf{U}\!B^{-1}$",align=E);

f.line(b, "$\textbf{b}$", position=0.95, align=(-1,0),draw_orientation=true,topair(R));
var bb = Rot*b/Rot;
f.line(bb, "$B\textbf{b}B^{-1}$", position=0.95, align=(-1,0),draw_orientation=true,topair(Rot*R/Rot) );

bool show_point_at_infinity = false;
if(show_point_at_infinity) {
  var N = wedge(e_0, Line(0,-1,0));
  f.point_at_infinity(N, "$\textbf{N}$", O=(3.5,1.5));
  f.point_at_infinity(Rot*N/Rot, "$B\textbf{N}B^{-1}$", align=(-1,0), O=(3.5,1.5));
}

\end{asy}
\end{subfloatenv}\hfill%
\begin{subfloatenv}{ }
\begin{asy}
import Figure2D;
Figure f = Figure(yaxis_name="$t$");
metric = Metric(Minkowski);
f.null_lines((-2,-1));

MV b = -Line(1,-2,-1/2);
b=b/norm(b);
var P = Point(1,1,-2);

f.line(b, "$\textbf{b}$", position=0.4, align=E);
f.point(P, "$\textbf{P}$");

MV a = Line(0,2,-1);
a = a/norm(a);
var N = wedge(e_0, a);
f.point_at_infinity(N, "$\textbf{e}_0\wedge\textbf{a}$", O=(-2,-1), align=(0,-1));

real tau = sqrt(3);
var Tran = exp(- 1/2 * tau * N);
f.line(Tran*b/Tran, "$T\textbf{b}T^{-1}$", position=0.3, align=E);
f.point(Tran*P/Tran, "$T\textbf{P}T^{-1}$", align=(0.8,1.1));

bool donull = false;
if(donull) {
MV a = Line(0,1,-1);
var N = wedge(e_0, a);
f.point_at_infinity(N, "$\textbf{e}_0\wedge\textbf{a}$", O=(-2,-1), align=(0,-1));

real tau = 1;
var Tran = exp(- 1/2 * tau * N);
f.line(Tran*b/Tran, "$T\textbf{b}T^{-1}$", position=0.3, align=E);
f.point(Tran*P/Tran, "$T\textbf{P}T^{-1}$", align=(0.8,1.1));
}

\end{asy}
\end{subfloatenv}
\caption{Rotation and translation in \M{2}}
\label{rotation and translation M2}
\end{figure}

The action \(BMB^{-1}\) of a spinor
\begin{equation}
B=e^{-\tfrac{1}{2}\phi_B\tb{S}}=
\cosh{\tfrac{\phi_B}{2}}-\tb{S}\sinh{\tfrac{\phi_B}{2}}
\end{equation}
on multivector \(M\) yields a rotation around 
a normalised clockwise point 
\(\tb{S}=-(\e_{12}+x_R\e_{20}+t_R\e_{01})\) by the angle \(\phi_B\).
It is illustrated in Figure~\ref{rotation and translation M2}(a),
where \(\tb{S}=-\e_{12}\) and \(\phi_B = \tfrac{3}{4}\).
A clockwise point  is used in the definition since  \(\phi\), which parametrises proper lines,
increases in the clockwise direction.
In Euclidean space \E{2}, the parametrising angle \(\alpha\) increases in the counterclockwise direction,
so a counterclockwise point is used in that case.

The rotation of a line parametrised by Equation~(\ref{line parameter phi})
corresponds to relativistic addition of velocities.
Accordingly, such rotation is usually called a boost.
A clockwise rotation (\(\phi_B>0\)) corresponds to a boost in the direction of increasing \(x\),
while a counterclockwise rotation (\(\phi_B<0\)) results in a boost in the direction of decreasing \(x\).
Let 
\begin{equation}
\tb{b}'=B\tb{b}B^{-1},
\end{equation} where \(\tb{b}\) and \(\tb{b}'\)
are parametrised by  \(\phi\) and \(\phi'\), respectively.
Then, 
\begin{equation}
\phi'=\phi+\phi_B.
\end{equation}
%If line \(\tb{a}\) represents a worldline of an object moving at velocity \(u=-\tanh{\phi}\)
%perceived by an observer who is at rest, then
%line \(\tb{a}'\) is a worldline of the same object, moving at velocity \(u'=-\tanh{\phi'}\),
%perceived by an observer at \(\tb{R}\) who is moving at velocity \(v=-\tanh{\phi_B}\).
Applying the addition formula for hyperbolic tangent,
\begin{equation}\label{tanh addition}
\tanh\phi'=\tanh{(\phi+\phi_B)}=\frac{\tanh{\phi}+\tanh{\phi_B}}{1+\tanh{\phi}\tanh{\phi_B}},
\end{equation}
gives the familiar expression for
relativistic addition of velocities
\begin{equation}
u'=\frac{u+v}{1+uv},
\end{equation}
where  \(u=\tanh{\phi}\), \(u'=\tanh{\phi'}\), and \(v=\tanh{\phi_B}\).
Note that \(u\) and \(v\) play distinct roles in this formula:
\(u\) is a velocity that is being acted upon and \(v\) is a velocity that characterises the action itself.
%In the context of special relativity,
%if \(\tb{b}\) is a worldline of an object perceived by an observer at \(\tb{R}\) who is at rest,
%then \(\tb{b}'\) is a worldline the same object perceived by another observer at \(\tb{R}\)
%moving at velocity \(-v=\tanh{\phi_B}\).

The rotation of a finite point corresponds to a Lorentz transformation.
In terms of coordinates in the standard basis, equation
\begin{equation}\tb{P}'=B\tb{P}B^{-1},
\end{equation} where 
\(\tb{P}=\e_{12}+x\e_{20}+t\e_{01}\) and
\(\tb{P}'=\e_{12}+x'\e_{20}+t'\e_{01}\),
gives
\begin{equation}
\left\{
\begin{aligned}
&x'=x\cosh{\phi_B}+t\sinh{\phi_B},\\
&t'=x\sinh{\phi_B}+t\cosh{\phi_B}.
\end{aligned}
\right.
\end{equation}
and, therefore,
\begin{equation}
\left\{
\begin{aligned}
&x'=\gamma\,(x+tv),\\
&t'=\gamma\,(t+xv),
\end{aligned}
\right.
\end{equation}
where \(\gamma=\cosh{\phi_B}=(1-\tanh^2\!\phi_B)^{-\tfrac{1}{2}}\) 
is the Lorentz factor corresponding to \(v=\tanh\phi_B\).
%In the context of special relativity,
%if \(\tb{P}\) is an event perceived by an observer  at rest,
%then \(\tb{P}'\) is the same event perceived by another observer
%moving at velocity \(-v=\tanh{\phi_B}\).

The points \(\tb{P}\), \(\tb{Q}\), \(\tb{R}\), \(\tb{U}\) shown in Figure~\ref{rotation and translation M2}(a) were chosen
such that the norm of their join with \(\tb{S}\) equals 2, e.g.\ \(\norm{\tb{P}\vee\tb{S}}=2\) (all points including \(\tb{S}\) are normalised).
Under the action of the spinor \(B=e^{-\tfrac{1}{2}\phi_B\tb{S}}\) with variable \(\phi_B\) ranging from \(-\infty\) to \(+\infty\),
the points \(\tb{P}\), \(\tb{Q}\), \(\tb{R}\), \(\tb{U}\)  follows trajectories in \M{2} shown in Figure~\ref{rotation and translation M2}(a)
with solid and dashed curves.
The points \(\tb{X}\) which lie on these curves satisfy \(\norm{\tb{X}\vee\tb{S}}=2\).
The two branches for which \(\tb{X}\vee\tb{S}\) is proper (solid curves) consist of points at the distance 2 from \(\tb{S}\).
I will refer to the shape defined by \(\norm{\tb{X}\vee\tb{S}}=r\), where \(\tb{X}\vee\tb{S}\) is proper, as a proper circle of radius \(r\) centred at \(\tb{S}\).
The proper circle consists of points at the distance \(r\) from the centre \(\tb{S}\).
It is convenient to refer to the shape defined by \(\norm{\tb{X}\vee\tb{S}}=r\), where \(\tb{X}\vee\tb{S}\) is \emph{improper},
as an improper circle of radius \(r\) centred at \(\tb{S}\).
However, in this case the term radius is used informally since the distance is undefined.
A point undergoing rotation stays on the same branch.
The line \(\tb{b}\) shown in  Figure~\ref{rotation and translation M2}(a) passes through \(\tb{R}\) and is tangent to the branch passing through \(\tb{R}\).
The rotated line \(\tb{b}'=B\tb{b}B^{-1}\) passes through \(B\tb{R}B^{-1}\) and is tangent to the same branch.

The action of a spinor
\begin{equation}
T=e^{-\tfrac{1}{2}\lambda\e_0\wedge\tb{a}}=1-\tfrac{\lambda}{2}\e_0\wedge\tb{a},
\end{equation}
where \(\tb{a}\) is a finite line and \(\lambda\in\R{}\),
on a point \(\tb{P}\) yields a translation of \(\tb{P}\) along the line perpendicular to \(\tb{a}\).
If \(\tb{a}\) is null, the translation is along  \(\tb{a}\), which is in accord with the observation that a null line is perpendicular to itself.
If \tb{a} is not null, it can be normalised.
If the normalised \(\tb{a}\) is improper, then the distance of translation is given by \(\lambda\) and \(\norm{\tb{P}\vee(T\tb{P}T^{-1})}=\lambda\) is satisfied
(\(\tb{P}\) is normalised and \(\lambda>0\)).
If the normalised \(\tb{a}\) is proper, then the distance of translation is undefined but \(\norm{\tb{P}\vee(T\tb{P}T^{-1})}=\lambda\) is still satisfied.
The same action applied to a line yields similar results.

In Euclidean case, the direction of translation is coincident with the top-down orientation vector of \(\tb{a}\) since it is perpendicular to \(\tb{a}\).
In \M{2}, the direction of translation is still perpendicular to \(\tb{a}\) but it does not necessarily coincide with the top-down orientation vector of \(\tb{a}\).
In fact, the direction of translation in \M{2} can be opposite to the direction of the top-down orientation vector.
For instance, a translation by \(\lambda>0\) along the \(t\)-axis in the direction of increasing \(t\) is induced by \(e^{-\tfrac{1}{2}\lambda\e_{20}}\)
while a translation along the \(x\)-axis in the direction of increasing \(x\) is induced by  \(e^{-\tfrac{1}{2}\lambda\e_{01}}\).
In the standard basis,  \(\tb{P}'=T\tb{P}T^{-1}\) gives
\begin{equation}
\left\{
\begin{aligned}
&x'=x + \lambda a,\\
&t'=t-\lambda h,
\end{aligned}
\right.
\label{Translation 1d}
\end{equation}
where \(\tb{a}=a\e_1+h\e_2\) and \(\e_0\wedge\tb{a}=-h\e_{20}+a\e_{01}\).
This is illustrated in Figure~\ref{rotation and translation M2}(b),
where \(\lambda=\sqrt{3}\),
\(\tb{a}=\tfrac{1}{\sqrt{3}}(2\e_1-\e_2)\) is normalised, and
 \(\e_0\wedge\tb{a}=\tfrac{1}{\sqrt{3}}(\e_{20}+2\e_{01})\) .

\section{Minkowski space \M{3} ( $2+1$ spacetime)}
In the standard notation for planes and points, I replace \(c\) with \(h\) and \(z\) with \(t\).
For
a plane \(\tb{a}=d\e_0+a\e_1+b\e_2+h\e_3\), 
a line \(\mb{\Lambda}=p_{10}\e_{10}+p_{20}\e_{20}+p_{30}\e_{30}+p_{23}\e_{23}+p_{31}\e_{31}+p_{12}\e_{12}\),
where \(p_{10}p_{23}+p_{20}p_{31}+p_{30}p_{12}=0\), and 
a point \(\tb{P}=w\e_{123}+x\e_{320}+y\e_{130}+t\e_{210}\),
I have
\begin{equation}
\tb{a}^2=a^2+b^2-h^2,\quad
\mb{\Lambda}^2=p_{23}^2+p_{31}^2-p_{12}^2,\quad 
\tb{P}^2=w^2,
\end{equation}
and
\begin{equation}
\norm{\tb{a}}=|a^2+b^2-h^2|^{\tfrac{1}{2}},\quad
\norm{\mb{\Lambda}}=|p_{23}^2+p_{31}^2-p_{12}^2|^{\tfrac{1}{2}},\quad
\norm{\tb{P}}=|w|.
\end{equation}
Points, lines, and the plane at infinity have zero norm.
A finite plane is proper if \(\tb{a}^2>0\), improper if  \(\tb{a}^2<0\), and null if \(\tb{a}^2=0\).
A finite line \(\mb{\Lambda}\) is proper if \(\mb{\Lambda}^2<0\) (note the reversal of the sign as compared with the definition of a proper plane),
improper if \(\mb{\Lambda}^2>0\), and null if  \(\mb{\Lambda}^2=0\).
Proper lines corresponds to worldlines of objects moving at less than the speed of light.
Null lines  passing through a given finite point form a pair of cones whose axes are parallel to the \(t\)-axis.
A null plane is tangent to a null cone.
A proper line \(\mb{\Lambda}\) is future-oriented if \(\mb{\Lambda}\cdot\e_{12}>0\) and past-oriented if \(\mb{\Lambda}\cdot\e_{12}<0\).
This is consistent with the terminology used in \M{2}, e.g., 
a proper line is future-oriented if its bottom-up orientation is facing in the direction of increasing \(t\).

\begin{figure}[h!]
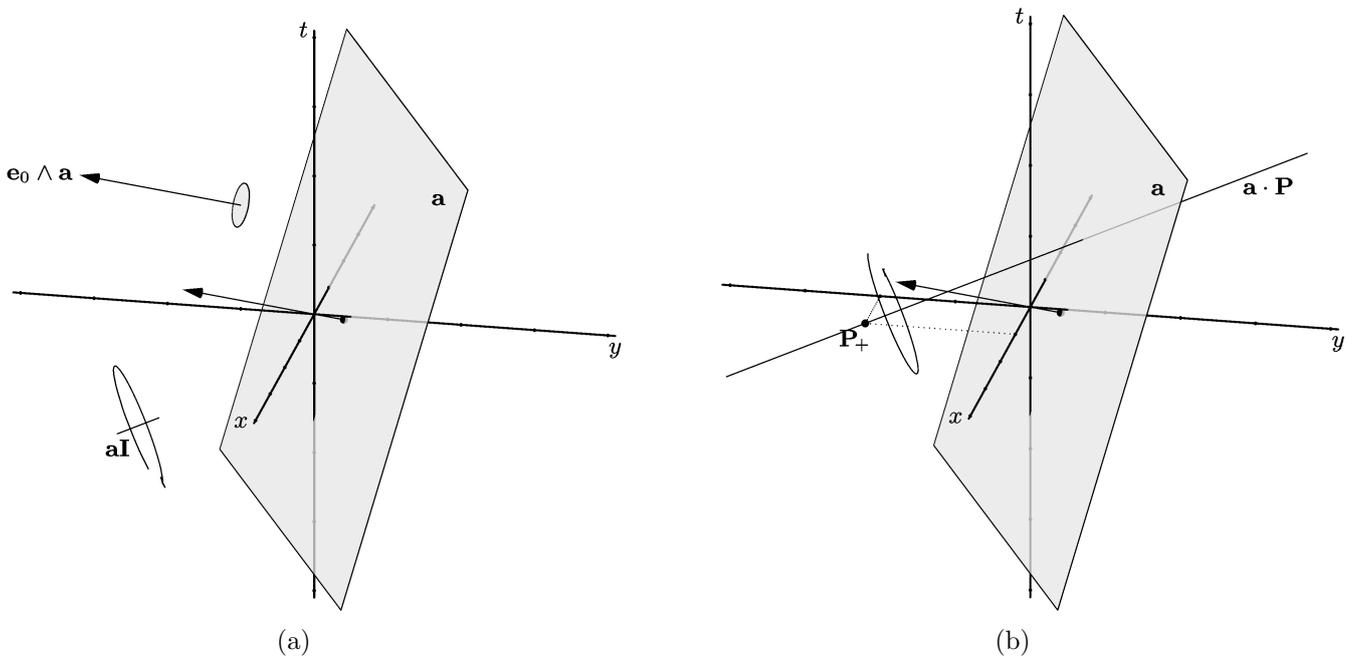
\hspace{-1cm}
\begin{subfloatenv}{ }
\begin{asy}
import Figure3D;
metric = Metric(Minkowski);
Figure f = Figure(metric);

var a = Plane(1,1,-2,1/1.5);
f.plane(a, "$\textbf{a}$", align=(1,-3,-1));

var L = wedge(e_0,a);
f.line_at_infinity(L, (0,-1,1.5), "$\textbf{e}_0\wedge\textbf{a}$", align=(0,-1,0));

var N = a*I;
f.point_at_infinity(N, (2,-2,-1), "$\textbf{a}\textbf{I}$", align=(0.5,-1.5,-2));

\end{asy}
\end{subfloatenv}\hfill%
\begin{subfloatenv}{ }
\begin{asy}
import Figure3D;
metric = Metric(Minkowski);
Figure f = Figure(metric);

var a = Plane(1,1,-2,1/1.5);
f.plane(a, "$\textbf{a}$", align=(1,-3,-1));

var P = Point(1,1,-2,0);
f.point(P, "$\textbf{P}$", align=(0,-0.5,-1));

var L = dot(a,P);
f.line(L, "$\textbf{a}\cdot\textbf{P}$", align=(0,-0.5,-2),position=0.05,orientation_shift=1.75);

//var Q = wedge(a,dot(a,P));
//f.stretch(P, Q, "$|\textbf{a}\vee\textbf{P}|$", label_angle=17,align=(-0.15,0.6), shift_distance=0.95);

\end{asy}
\end{subfloatenv}
\caption{Basic properties of points and planes  in \M{3}}
\label{basic  M3}
\end{figure}

The line at infinity \(\e_0\wedge\tb{a}\) provides the top-down orientation of a plane \(\tb{a}\).
The orientation is depicted in Figure~\ref{basic  M3}(a) by an arrow attached to the central point of the plane.
It only indicates the orientation of the relevant stack of planes.
The direction perpendicular to the plane is given by the polar point \(\tb{a}\I\).
The same direction is given by the line \(\tb{a}\cdot\tb{P}\), which passes through \(\tb{P}\)
as shown in Figure~\ref{basic  M3}(b).
If \(\tb{a}\) and \(\tb{P}\) are normalised and \(\tb{a}\) is improper, which implies \(\tb{a}\cdot\tb{P}\) is proper, then \(|\tb{a}\vee\tb{P}|\) gives the distance 
between \(\tb{a}\) and \(\tb{P}\).
Note that \(\tb{a}\) shown in Figure~\ref{basic  M3}(b) is proper, hence the distance between \(\tb{a}\) and \(\tb{P}\) is undefined.

\begin{figure}[t!]
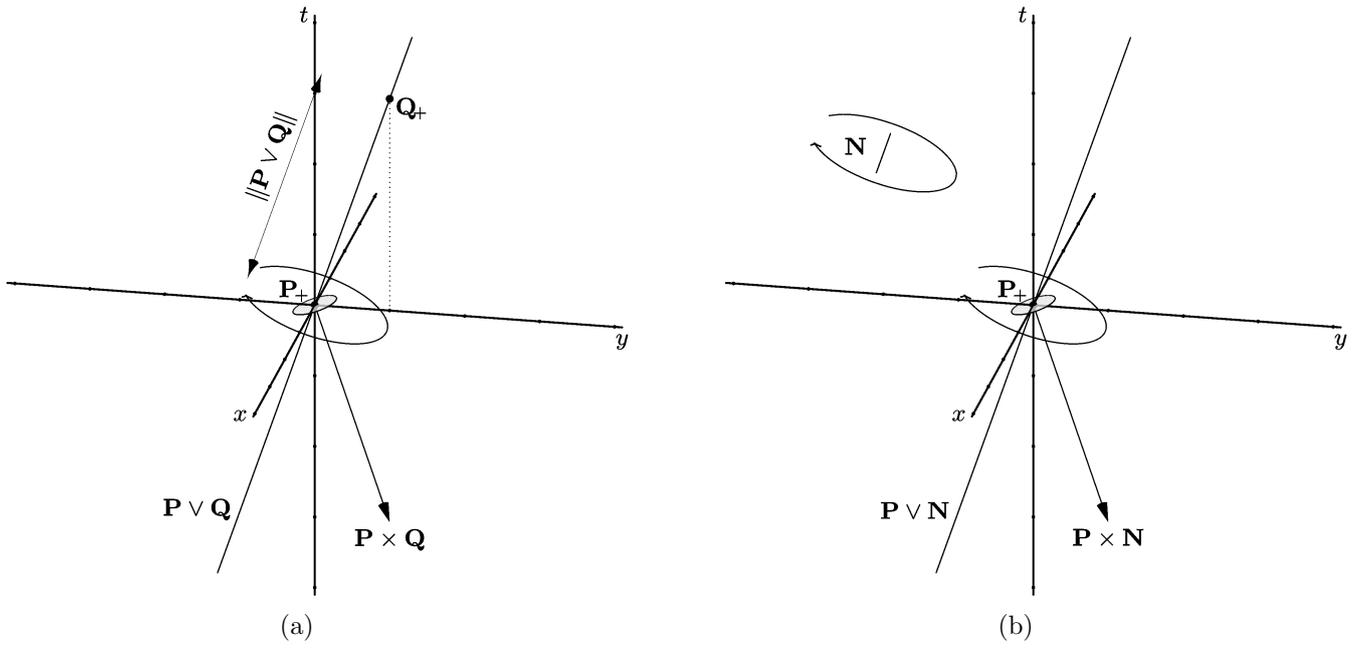
\hspace{-1cm}
\begin{subfloatenv}{ }
\begin{asy}
import Figure3D;
metric = Metric(Minkowski);
Figure f = Figure(metric);

var P = Point(1,0,0,0);
var Q = Point(1,0,1,3);
f.point(P, "$\textbf{P}$", align=(0,-1,0.75));
f.point(Q, "$\textbf{Q}$", align=(0,1,-0.5));

var L = join(P,Q);
f.line(L, "$\textbf{P}\vee\textbf{Q}$", position=0.9,align=(0,-1,0.5));

var L0 = cross(P,Q);
f.line_at_infinity(L0, (0,0,0), "$\textbf{P}\times\textbf{Q}$", align=(0,0,-1));

f.stretch(P, Q, "$\lVert\textbf{P}\vee\textbf{Q}\rVert$", label_angle=70,align=(-0.85,0.6), shift_distance=0.95);

\end{asy}
\end{subfloatenv}\hfill%
\begin{subfloatenv}{ }
\begin{asy}
import Figure3D;
metric = Metric(Minkowski);
Figure f = Figure(metric);

var P = Point(1,0,0,0);
var N = wedge(e_0, -Line(0,0,0,0,1,3));
f.point(P, "$\textbf{P}$", align=(0,-1,0.75));
f.point_at_infinity(N, (0,-2,2), "$\textbf{N}$", align=(0,-2.5,0.5));

var L = join(P,N);
f.line(L, "$\textbf{P}\vee\textbf{N}$", position=0.9,align=(0,-1,0.5));

var L0 = cross(P,N);
f.line_at_infinity(L0, (0,0,0), "$\textbf{P}\times\textbf{N}$", align=(0,0,-1));

\end{asy}
\end{subfloatenv}
\caption{The join of two points in \M{3}}
\label{basic  M3 2}
\end{figure}

\begin{figure}[t!]
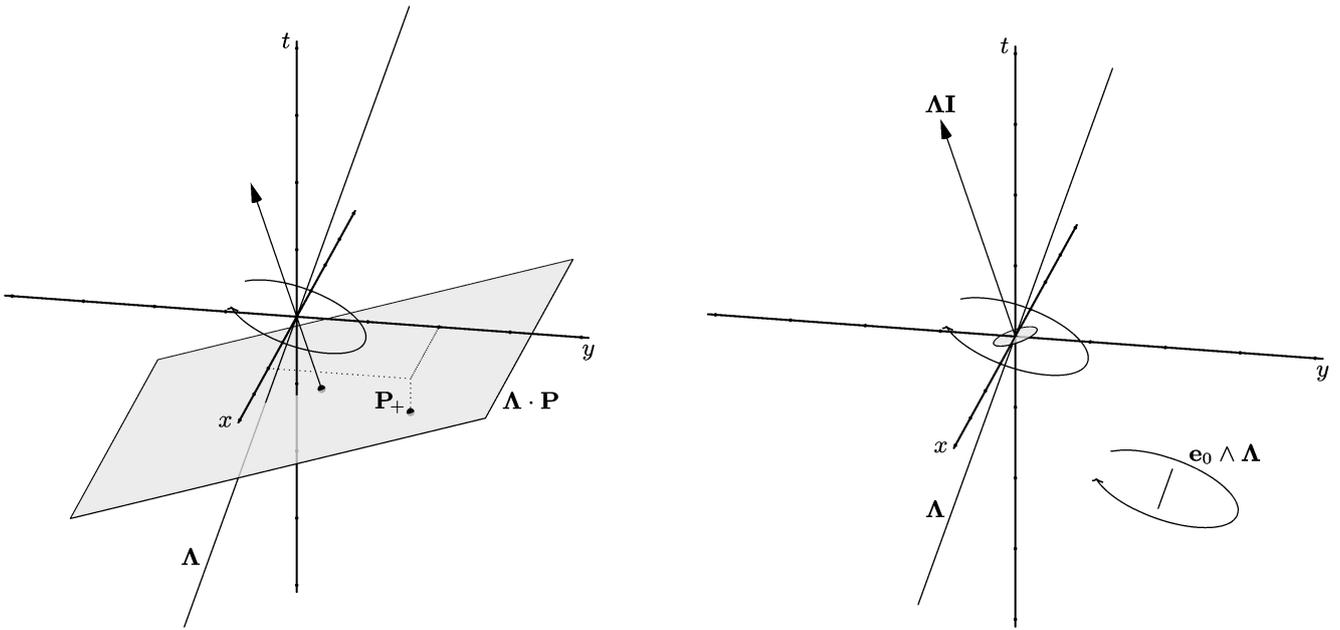
\hspace{-1cm}
\begin{subfloatenv}{}
\begin{asy}
import Figure3D;
metric = Metric(Minkowski);
Figure f = Figure(metric);

var R = Point(1,0,0,0);
var Q = Point(1,0,1,3);
var L = join(R,Q);

var P = Point(1,2,2,-1/2);

f.line(L, size=5,"$\mathbf{\Lambda}$", position=0.9, align=(0,-1,0.5));
f.point(P, "$\textbf{P}$",align=(0,-1,0.5));

f.plane(dot(L,P),"$\mathbf{\Lambda}\cdot\textbf{P}$", align=(0,2.5,2));

\end{asy}
\end{subfloatenv}\hfill%
\begin{subfloatenv}{}
\begin{asy}
import Figure3D;
metric = Metric(Minkowski);
Figure f = Figure(metric);

var R = Point(1,0,0,0);
var Q = Point(1,0,1,3);
var L = join(R,Q);

f.line(L, "$\mathbf{\Lambda}$", position=0.85, align=(0,-0.75,1));

f.point_at_infinity(wedge(e_0,L),"$\textbf{e}_0\wedge\mathbf{\Lambda}$", centre=(0,2,-2), align=(0,3.5,4));

f.line_at_infinity(L*I,"$\mathbf{\Lambda}\textbf{I}$", centre=(0,0,0),align=(0,0,1));

\end{asy}
\end{subfloatenv}
\caption{Points and lines in \M{3}}
\label{basic lines and points in M3}
\end{figure}

The commutator \(\tb{P}\times\tb{Q}\) is a line at infinity whose stack consists of the planes perpendicular 
to  \(\tb{P}\vee\tb{Q}\) (see Figure~\ref{basic  M3 2}(a)).
The same applies to the commutator \(\tb{P}\times\tb{N}\), where \(\tb{N}\) is a point at infinity.
The distance \(r\)  between normalised points \(\tb{P}\) and \(\tb{Q}\) connected by a proper line
is given by 
\begin{equation}
r=\norm{\tb{P}\vee\tb{Q}},
\end{equation}
where \(\tb{P}\vee\tb{Q}\) is a line passing through  \(\tb{P}\) and \(\tb{Q}\).
In the standard basis,  \(\tb{P}=\e_{123}+x_P\e_{320}+y_P\e_{130} +t_P\e_{210}\) and 
\(\tb{Q}=\e_{123}+x_Q\e_{320}+y_Q\e_{130} +t_Q\e_{210}\), which gives
\(r=\sqrt{(t_P-t_Q)^2-(x_P-x_Q)^2-(y_P-y_Q)^2}\).
If \(\tb{P}\vee\tb{Q}\) is null or improper, the distance is undefined.
However, formally \(r=0\) if the points are connected by a null line.
As in \E{3}, the join of four points is a scalar and, if the points are normalised,
\(\tfrac{1}{3!}|\tb{P}\vee\tb{Q}\vee\tb{R}\vee\tb{S}|\) gives the volume of a 3-simplex 
defined by the points. Its value is identical to that in \E{3}.

The dot product \(\mb{\Lambda}\cdot\tb{P}\) of a line \(\mb{\Lambda}\) and a point \(\tb{P}\) 
is a plane that passes through \(\tb{P}\) and is perpendicular to \(\mb{\Lambda}\),
as illustrated in Figure~\ref{basic lines and points in M3}(a).
The commutator \(\mb{\Lambda}\times\tb{P}\) is a point at infinity that is perpendicular to the plane 
\(\mb{\Lambda}\vee\tb{P}\).
The distance between \(\mb{\Lambda}\) and \(\tb{P}\)
is defined if \(\mb{\Lambda}_P\) is proper, 
where \(\mb{\Lambda}_P=(\mb{\Lambda}\cdot\tb{P})\wedge(\mb{\Lambda}\vee\tb{P})\) passes through \(\tb{P}\) and is perpendicular to \(\mb{\Lambda}\).
When defined, the distance is given by \(\norm{\mb{\Lambda}\vee\tb{P}}\) if \(\mb{\Lambda}\) and \(\tb{P}\) are normalised.
The dot product \(\tb{a}\cdot\mb{\Lambda}\) of a plane \(\tb{a}\) and a line \(\mb{\Lambda}\)
is a plane that passes through \(\mb{\Lambda}\) and is perpendicular to \(\tb{a}\).
The line at infinity \(\mb{\Lambda}\I\) is perpendicular to \(\mb{\Lambda}\), i.e.\ 
its stack consists of the planes perpendicular to \(\mb{\Lambda}\).

Essentially, the above results are identical to those in \E{3}, except that the notion of perpendicularity
is different in \M{3}.
The other notable distinction is that there are finite lines and planes in \M{3} which are not invertible, the null objects.

\begin{figure}[t!]
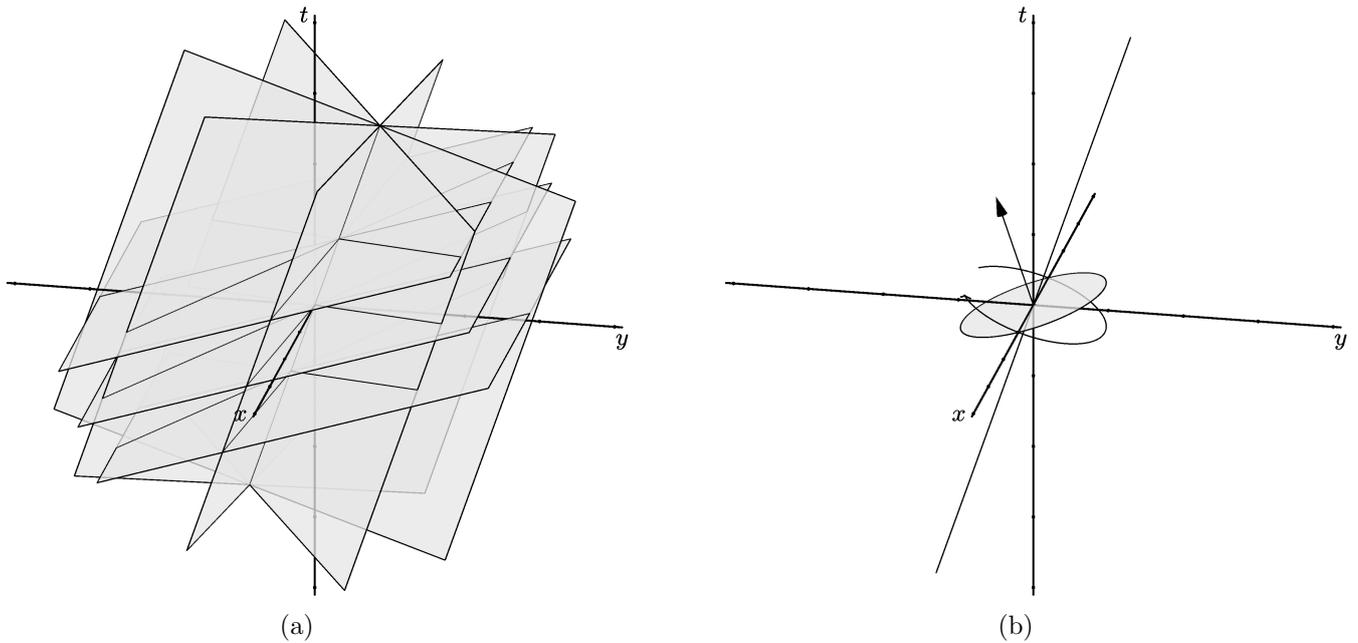
\hspace{-1cm}
\begin{subfloatenv}{ }
\begin{asy}
import Figure3D;
metric = Metric(Minkowski);
Figure f = Figure(metric);

var R = Point(1,0,0,0);
var Q = Point(1,0,1,3);
var L = join(R,Q);

var C = L + I*L/2;

var L = (1 - 0.5*join(C,C)/dot(C,C)*I)*C;
var Lp = 0.5*join(C,C)/dot(C,C)*I*C;

L = L/norm(L);
triple n = (L.p23,L.p31,L.p12);
var a = join(L,Point(1,0,0,1));
int N = 4;

metric = Metric(Euclidean);
for(int i: sequence(N)) { real alpha = i/N*pi; var Rot = cos(alpha/2) - sin(alpha/2)*L; f.plane(Rot*a/Rot,u=n,draw_orientation=false,size=2.75);}
metric = Metric(Minkowski);

var b = join(O,L*I);
int M = 1;
metric = Metric(Euclidean);
for(int i: sequence(-M,M)) { real alpha = i/M*pi; var Tr = 1 - i/2*Lp/1.95; write(Tr); f.plane(Tr*b/Tr,draw_orientation=false,size=2.75);}

\end{asy}
\end{subfloatenv}%\hfill%
\begin{subfloatenv}{ }
\begin{asy}
import Figure3D;
metric = Metric(Minkowski);
Figure f = Figure(metric);

var R = Point(1,0,0,0);
var Q = Point(1,0,1,3);
var L = join(R,Q);

var cL = L + I*L/2;

f.compound_line(cL);

\end{asy}
\end{subfloatenv}
\caption{Visualising a non-simple bivector in \M{3}}
\label{non-simple bivector in M3}
\end{figure}

As in \E{3}, any non-simple bivector \(\mb{\Lambda}\) in \M{3} can be decomposed into a line \(\mb{\Lambda}_0\)
that passes through the origin and a line at infinity \(\mb{\Lambda}_\infty=\e_0\wedge\tb{a}\),
where  \(\tb{a}\) can be assumed to pass through the origin as well.
If \(\mb{\Lambda}_0\) is not null, then \(\mb{\Lambda}\cdot\mb{\Lambda}\ne0\) and the same reasoning that was used for 
non-simple bivectors in \E{3} applies in \M{3}.
So, I can write \(\mb{\Lambda}=\mb{\Lambda}_1 + \mb{\Lambda}_2\), where 
\(\mb{\Lambda}_1\) is the finite axis of the bivector
and  \(\mb{\Lambda}_2\) is the axis at infinity, which is perpendicular to the finite axis
in the sense that the stack of \(\mb{\Lambda}_2\) consists of the planes perpendicular to \(\mb{\Lambda}_1\).
As in Euclidean space, I have
\begin{equation}
\mb{\Lambda}_1=(1-a\I)\mb{\Lambda},\quad
\mb{\Lambda}_2=a\I\mb{\Lambda},\quad
\textrm{and } 
a=\frac{\mb{\Lambda}\vee\mb{\Lambda}}{2\mb{\Lambda}\cdot\mb{\Lambda}}.
\end{equation}
An example is shown in Figure~\ref{non-simple bivector in M3}.
If \(\mb{\Lambda}_0\) is null, this convenient decomposition is no longer valid since \(\mb{\Lambda}\cdot\mb{\Lambda}=0\).

The commutator of two parallel lines is a line at infinity whose stack is perpendicular to the plane passing
through both parallel lines.
The commutator of two lines intersecting at a finite point is a finite line that passes through the point
of the intersection and is perpendicular to the plane passing through the intersecting lines.
However, if the plane passing through the intersecting lines is null, the commutator is a null line which
lies in this plane. 
In general, the commutator of two lines is a non-simple bivector.
If the commutator is null, is cannot be decomposed into orthogonal axes.

\begin{figure}[t!]
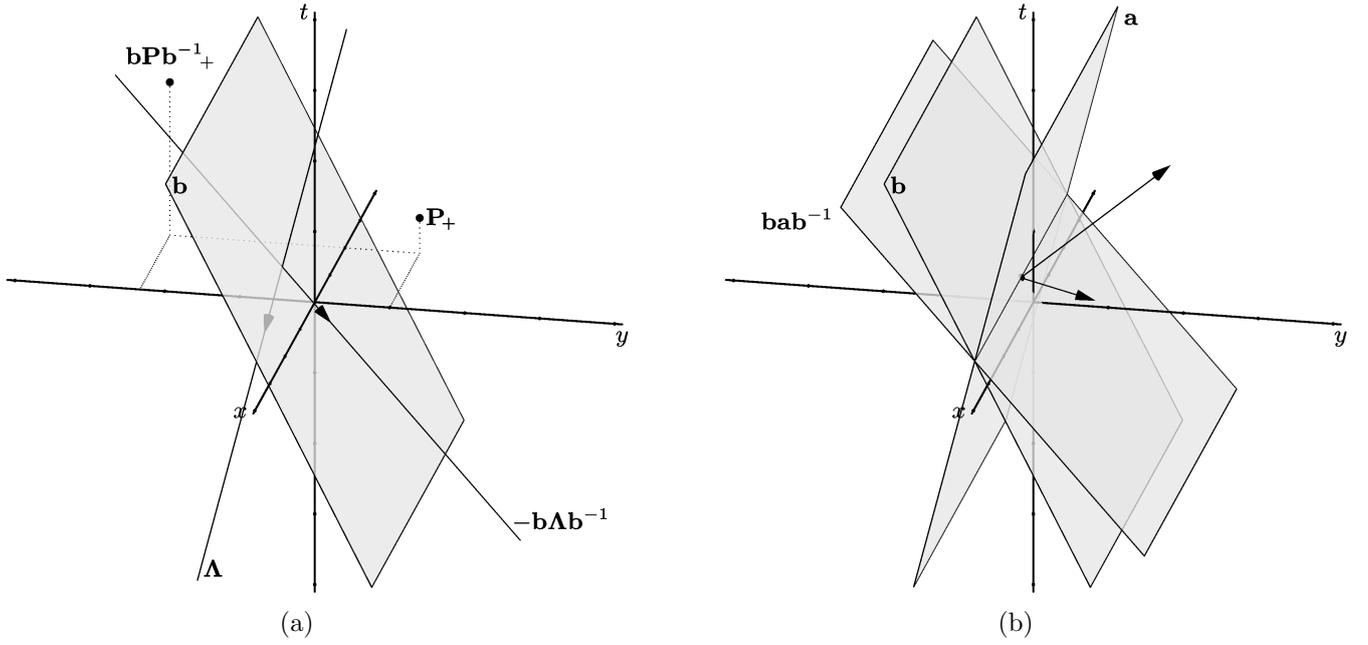

\hspace{-1cm}
\begin{subfloatenv}{ }
\begin{asy}
import Figure3D;
metric = Metric(Minkowski);
Figure f = Figure(metric);

var a = Plane(0,0,2,1);
f.plane(a,"$\textbf{b}$", draw_orientation=false, align=(0,1,0),label_position=3);

var P = Point(1,-2,1,1/2);
f.point(P,"$\textbf{P}$", draw_helper_lines=true);
f.point(a*P/a, "$\textbf{b}\textbf{P}\textbf{b}^{-1}$", align=(0,0,2));

var L = Line(2,2,-1/2,0,1,4)/8;

triple A = totriple(wedge(a,L));
f.line(L,"$\mathbf{\Lambda}$",position=0.05,align=(0,0.5,-1),orientation_theta=-10,draw_orientation=false,draw_bottom_up_orientation=true,O=A);
f.line(-a*L/a, "$-\textbf{b}\mathbf{\Lambda}\textbf{b}^{-1}$",position=0.05,align=(0,2,0),orientation_theta=-10,draw_orientation=false,draw_bottom_up_orientation=true,O=A);

//var N = wedge(e_0,L);
//f.point_at_infinity(N,(0,-3,-2.5),"$\textbf{N}$",align=(0,4,-2),orientation_theta=-10);
//f.point_at_infinity(-a*N/a,(0,-3,-1.5),"$-\textbf{b}\textbf{N}\textbf{b}^{-1}$",align=(0,-1,4));

\end{asy}
\end{subfloatenv}%\hfill%
\begin{subfloatenv}{ }
\begin{asy}
import Figure3D;
metric = Metric(Minkowski);
Figure f = Figure(metric);

var a = Plane(0,0,2,1);
f.plane(a,"$\textbf{b}$", draw_orientation=false, align=(0,1,0),label_position=3);

//var K = Line(0,-4,1,0,0,0)/4;
//f.line_at_infinity(K,(0,3,3.5),"$\mathbf{\Phi}$",align=(0,-0.5,-1));
//f.line_at_infinity(a*K/a,(0,3,3.5), "$\textbf{b}\mathbf{\Phi}\textbf{b}^{-1}$",align=(0,2,-1));

var b = Plane(1,0,4,-1)/4;

MV C = wedge(wedge(a,b),e_1);
triple c = totriple(C);
f.plane(b,u=(0,1,0),"$\textbf{a}$",draw_orientation=true,O=c,align=(0,1,-1));
f.plane(a*b/a,u=(0,1,0),"$\textbf{b}\textbf{a}\textbf{b}^{-1}$",label_position=1,align=(0,-1.5,-1),draw_orientation=true,O=c);

\end{asy}
\end{subfloatenv}
\caption{Bottom-up reflection in a plane in \M{3}}
\label{reflection in a plane in M3}
\end{figure}

Projection, rejection, reflection, and scaling with respect to an invertible blade \(A_k\) are defined in the same way as in \E{3}.
Moreover, if \(A_k\) is a point, then the result of these transformations is the same as in \E{3}.
As in \E{3}, geometric objects in \M{3} are projected and rejected along the perpendicular direction.
So, these transformations can be readily understood by adopting the appropriate view of perpendicularity.
As in \E{3}, there are two alternative ways to project and reject a line with respect to another line.
Bottom-up reflections in a plane are illustrated in Figure~\ref{reflection in a plane in M3}.
Note that in Figure~\ref{reflection in a plane in M3}(b) the top-down orientation of the planes is shown, even though
the reflection is bottom-up.

A normalised proper line passing through the origin can be parametrised by the angles \(\alpha\) and \(\phi\) as follows
\begin{equation}
\mb{\Lambda}=
\e_{23}\cos{\alpha}\sinh{\phi}+\e_{31}\sin{\alpha}\sinh{\phi}+\e_{12}\cosh{\phi}.
\label{worldline in M3}
\end{equation}
It can be interpreted as the worldline of an object moving in the direction \(\alpha\), i.e.\ at the angle \(\alpha\) to the \(x\)-axis,
at the speed \(u=\tanh{\phi}\) if \(\phi>0\) (or in the direction \(-\alpha\) if \(\phi<0\)).
Note that the bottom-up orientation of \(\mb{\Lambda}\) defined by (\ref{worldline in M3}) is facing in the direction of decreasing \(t\),
i.e.\ \(\mb{\Lambda}\) is past-oriented (\(\mb{\Lambda}\cdot\e_{12}=-\cosh{\phi}<0\)).
If the opposite bottom-up orientation is required, it is provided by \(-\mb{\Lambda}\), which is future-oriented.
Whether \(\mb{\Lambda}\) or \(-\mb{\Lambda}\) is used has no effect on the formulas derived below.

The usual definition of the Spin group applies in \M{3}.
Namely, the Spin group consists of multivectors that can be written as the product of an even number
of normalised proper planes. Spinors are even multivectors and each spinor satisfies \(S\reverse{S}=1\).
The Lie algebra of the Spin group consists of bivectors, with the commutator used for multiplication.
Any proper motion in \M{3} can be generated by the action of  \(e^A\) for some bivector \(A\).
The action of
\begin{equation}
B=e^{-\tfrac{1}{2}\phi_B\mb{\Omega}_B},
\textrm{ where}\quad
\mb{\Omega}_B=-\e_{23}\sin{\alpha_B}+\e_{31}\cos{\alpha_B},
\label{booster in M3}
\end{equation}
on a proper line \(\mb{\Lambda}\) corresponds to a boost in the direction \(\alpha_B\) if \(\phi_B>0\).
Note that \(\mb{\Omega}_B\) is not a worldline.
Rather, it is a normalised bivector that enables the relevant rotation and a bivector happens to be a line in \M{3}.
In the standard basis, the equation
\begin{equation}
\mb{\Lambda}'=B\mb{\Lambda}B^{-1},
\label{spinor velocity transform in M3}
\end{equation} 
where \(\mb{\Lambda}'\) is parametrised with \(\alpha'\) and \(\phi'\),
%\begin{equation}
%\mb{\Lambda}'=\e_{23}\cos{\alpha'}\sinh{\phi'}+\e_{31}\sin{\alpha'}\sinh{\phi'}+\e_{12}\cosh{\phi'},
%\label{rotated worldline in M3}
%\end{equation} 
%is a worldline of the same object perceived by
%another observer moving in direction \(\alpha_R\) at speed \(\tanh{\phi_R}\).
%The worldline of the moving observer is given by
%\begin{equation}
%\mb{\Lambda}_R=
%\e_{23}\cos{\alpha_R}\sinh{\phi_R}+\e_{31}\sin{\alpha_R}\sinh{\phi_R}+\e_{12}\cosh{\phi_R}.
%\end{equation}
corresponds to 
\begin{equation*}
\!\!\!\left\{\!\!
\begin{aligned}
&\cos{\alpha'}\sinh{\phi'}=
\cos{\alpha}\sinh{\phi}
+\cos{\alpha_B}\cosh{\phi}\sinh{\phi_B}
+\cos{\alpha_B}\cos{(\alpha-\alpha_B)}\sinh{\phi}\,(\cosh{\phi_B}-1),\\
&\sin{\alpha'}\sinh{\phi'}=
\sin{\alpha}\sinh{\phi}
+\sin{\alpha_B}\cosh{\phi}\sinh{\phi_B}
+\sin{\alpha_B}\cos{(\alpha-\alpha_B)}\sinh{\phi}\,(\cosh{\phi_B}-1),\\
&\cosh{\phi'}=
\cosh{\phi}\cosh{\phi_B}+
\cos{(\alpha-\alpha_B)}\sinh{\phi}\sinh{\phi_B}.
\end{aligned}
\!\right.
\end{equation*}
Dividing the first two equations by the third and performing elementary substitutions
 yields 
%the expressions for relativistic transformation of velocity in the standard basis:
\[
\cos{\alpha'}\tanh{\phi'}
=
\frac{\dfrac{\cos{\alpha}\tanh{\phi}}{\cosh{\phi_B}}
+\cos{\alpha_B}\tanh{\phi_B}
+\dfrac{\cos{\alpha_B}\tanh{\phi_B}\cos{(\alpha-\alpha_B)}\tanh{\phi}\tanh{\phi_B}
\cosh{\phi_B}}{\cosh{\phi_B}+1}}
{1+
\cos{(\alpha-\alpha_B)}\tanh{\phi}\tanh{\phi_B}},
\]
\[\sin{\alpha'}\tanh{\phi'}
=
\frac{\dfrac{\sin{\alpha}\tanh{\phi}}{\cosh{\phi_B}}
+\sin{\alpha_B}\tanh{\phi_B}
+\dfrac{\sin{\alpha_B}\tanh{\phi_B}\cos{(\alpha-\alpha_B)}\tanh{\phi}\tanh{\phi_B}
\cosh{\phi_B}}{\cosh{\phi_B}+1}}
{1+
\cos{(\alpha-\alpha_B)}\tanh{\phi}\tanh{\phi_B}}.
\]
Let \((u_x,u_y)\) and \((u_x',u_y')\) be velocities associated with the worldlines 
 \(\mb{\Lambda}\) and \(\mb{\Lambda}'\), respectively, 
and \((v_x,v_y)\) be a velocity that characterises the boost.
Then 
\[
\begin{aligned}
&(u_x,u_y)=(\cos{\alpha}\tanh{\phi},\,\sin{\alpha}\tanh{\phi}),\\
&(u_x',u_y')=(\cos{\alpha'}\tanh{\phi'},\,\sin{\alpha'}\tanh{\phi'}),\\
&(v_x,v_y)=(\cos{\alpha_B}\tanh{\phi_B},\,\sin{\alpha_B}\tanh{\phi_B}),
\end{aligned}
\]
and, therefore,
\begin{equation}
\left\{
\begin{aligned}
&u'_x
=
\frac{
{\gamma}^{-1}u_x
+v_x
+v_x(u_xv_x+u_yv_y)\gamma(\gamma+1)^{-1}
}
{1+(u_xv_x+u_yv_y)},\\
&u'_y
=
\frac{
{\gamma}^{-1}u_y
+v_y
+v_y(u_xv_x+u_yv_y)\gamma(\gamma+1)^{-1}
}
{1+(u_xv_x+u_yv_y)},\\
\end{aligned}
\right.
\label{velocity transform in M3 components}
\end{equation}
where  \(\gamma=\cosh{\phi_B}\) is the Lorentz factor corresponding to \(v=\sqrt{v_x^2+v_y^2}\).

The action of the spinor \(B\) defined by (\ref{booster in M3})  on a finite point \(\tb{P}=\e_{123}+x\e_{320}+y\e_{130}+t\e_{210}\)
corresponds to a Lorentz transformation.
Substituting \(B\) into 
\begin{equation}
\tb{P}'=B\tb{P}B^{-1},
\end{equation}
where \(\tb{P}'=\e_{123}+x'\e_{320}+y'\e_{130}+t'\e_{210}\),
gives
\begin{equation*} 
\left\{
\begin{aligned}
&x'=x
+t\cos\alpha_B\sinh\phi_B
+x\cos^2\!\alpha_B\,(\cosh{\phi_B}-1)
+y\sin\alpha_B\cos\alpha_B\,(\cosh{\phi_B}-1),\\
&y'=y
+t\sin\alpha_B\sinh\phi_B
+x\sin\alpha_B\cos\alpha_B\,(\cosh{\phi_B}-1)
+y\sin^2\!\alpha_B\,(\cosh{\phi_B}-1), \\
&t'=t\cosh\phi_B + x\cos\alpha_B\sinh\phi_B+y\sin\alpha_B\sinh\phi_B,
\end{aligned}
\right.
\end{equation*}
and
\begin{equation} 
\left\{
\begin{aligned}
&x'=x
+\gamma tv_x
+(\gamma-1)\frac{(xv_x+yv_y)v_x}{v_x^2+v_y^2},\\
&y'=y
+\gamma tv_y
+(\gamma-1)\frac{(xv_x+yv_y)v_y}{v_x^2+v_y^2}, \\
&t'=\gamma(t + xv_x+yv_y).
\end{aligned}
\right.
\end{equation}

The bivector \(\mb{\Omega}_B\) defined by (\ref{booster in M3}) dually represents a line in \M{3}
that lies in \(\e_3\), the \(xy\)-plane, and passes through the origin.
A rotation around such a line corresponds to a boost in special relativity.
However, the boost defined by (\ref{spinor velocity transform in M3}) is not restricted 
to rotations around lines passing through the origin.
In fact, a rotation around any line \(\mb{\Omega}_B\) 
parallel to  the \(xy\)-plane, which implies \((\e_3\wedge\mb{\Omega}_B)^2=0\), corresponds to a boost.
Moreover, Equation (\ref{spinor velocity transform in M3}) can be applied to
any line \(\mb{\Lambda}\), even if it does not pass through the origin.

The action of 
\begin{equation}
R=e^{-\tfrac{1}{2}\alpha\mb{\Omega}_E},
\textrm{ where}\quad
\mb{\Omega}_E=\e_{12},
\label{Euclidean rotation in M3}
\end{equation}
corresponds to a Euclidean rotation around the line \(\e_{12}\), the \(t\)-axis, by the angle \(\alpha\).
For any line \(\mb{\Omega}_E\) perpendicular to the \(xy\)-plane, which implies \(\e_3\cdot\mb{\Omega}_E=0\),
the corresponding spinor \(R\) yields a Euclidean rotation around \(\mb{\Omega}_E\).

In general, the action of 
\begin{equation}
S=e^{-\tfrac{1}{2}\theta\mb{\Omega}},
\label{boost and Euclidean rotation in M3}
\end{equation}
where \(\mb{\Omega}\) is an arbitrary normalised or null line and \(\theta\in\R{}\),
corresponds to a general rotation in \M{3} around the line \(\mb{\Omega}\) by the angle \(\theta\).
It is neither a boost nor a Euclidean rotation.
Likewise,  \(\theta\) is neither purely Euclidean nor purely pseudo-Euclidean.
Instead, it encompasses two characteristics, speed and Euclidean angle, in a single quantity.
The character of the action of the spinor \(S\) depends on whether \(\mb{\Omega}\) is proper, null, or improper.
If \(\mb{\Omega}\) is proper, the Euclidean rotation dominates and the action of \(S\) is elliptic.
If \(\mb{\Omega}\) is improper, the boost dominates and the action is hyperbolic.
If \(\mb{\Omega}\) is null, then the action is parabolic.

\begin{figure}[t]
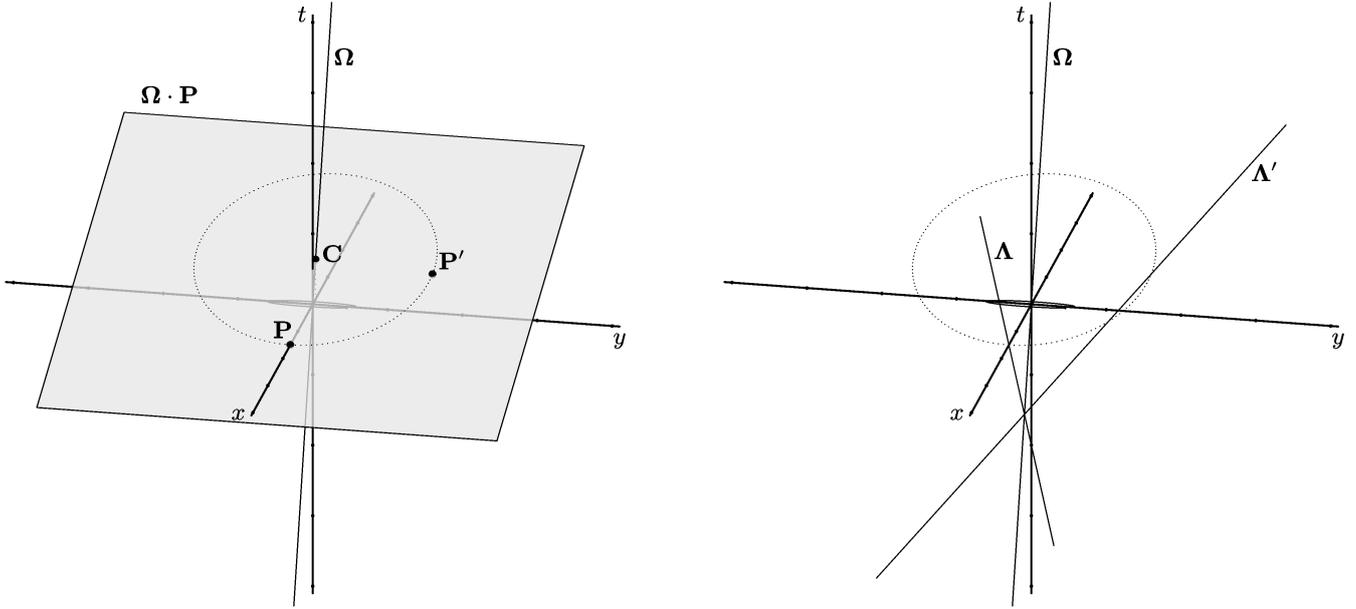
\hspace{-1cm}
\begin{subfloatenv}{}
\begin{asy}
import Figure3D;
metric = Metric(Minkowski);
Figure f = Figure(metric);

MV L = Line(0,0,0,-1,0,3);
L = L/norm(L);

MV P = Point(1,1.5,0,0);

f.line(L, size=4,"$\mathbf{\Omega}$", position=0.9, align=(0,1,0.5));
f.point(P, "$\textbf{P}$",align=(0,-0.5,1),draw_orientation=false);

f.plane(dot(L,P),"$\mathbf{\Omega}\cdot\textbf{P}$", align=(0,2.5,2),draw_orientation=false);

//int n = 10; for(int i: sequence(n)) { real theta = 2*pi*i/n; write(theta); MV S = exp(-1/2*theta*L); f.point(S*P/S,string(i), draw_helper_lines=false, draw_orientation=false); }
real theta = 1.5;
MV S = exp(-1/2*theta*L); 
f.point(S*P/S,"$\textbf{P}'$", draw_helper_lines=false, draw_orientation=false,align=(0,1,1));

MV C = wedge(L,dot(L,P));
f.point(C, "$\textbf{C}$",align=(1,1.5,1),draw_orientation=false);

//triple f(real theta) { MV S = exp(-1/2*theta*L); return totriple(S*P/S); }
triple f(real theta) { MV S = cos(-1/2*theta)+L*sin(-1/2*theta); return totriple(S*P/S); }
real x(real theta) { return f(theta).x; }
real y(real theta) { return f(theta).y; }
real z(real theta) { return f(theta).z; }
path3 p = graph(x, y, z, -pi, pi, operator ..);
draw(p,dotted);

//triple c = totriple(C);
//real xc = c.x; real yc = c.y; real tc = c.z;
//real r = norm(join(P,C));
//real g(pair p) { real x = p.x; real y = p.y; real d2 = (x-xc)^2+(y-yc)^2-r^2; if (d2>0) {return tc+sqrt(d2);} else {return 0;} }
//surface s = surface(g,(-2.5,-2.5),(2.5,2.5),nx=20);
//draw(s,lightgray,meshpen=black+thick(),nolight,render(merge=true));

\end{asy}
\end{subfloatenv}\hfill%
\begin{subfloatenv}{}
\begin{asy}
import Figure3D;
metric = Metric(Minkowski);
Figure f = Figure(metric);

MV L = Line(0,0,0,-1,0,3);
L = L/norm(L);
f.line(L, size=4,"$\mathbf{\Omega}$", position=0.9, align=(0,1,0.5));

var M = join(Point(1,1.5,0,0), Point(1,0,0,-2))/12;

f.line(M,"$\mathbf{\Lambda}$",position=0.9,draw_orientation=false);

real theta = 1.5;
MV S = exp(-1/2*theta*L); 
f.line(S*M/S,"$\mathbf{\Lambda}'$",position=0.9,draw_orientation=false);

MV P = Point(1,1.5,0,0);
//triple f(real theta) { MV S = exp(-1/2*theta*L); return totriple(S*P/S); }
triple f(real theta) { MV S = cos(-1/2*theta)+L*sin(-1/2*theta); return totriple(S*P/S); }
real x(real theta) { return f(theta).x; }
real y(real theta) { return f(theta).y; }
real z(real theta) { return f(theta).z; }
path3 p = graph(x, y, z, -pi, pi, operator ..);
draw(p,dotted);

\end{asy}
\end{subfloatenv}
\caption{Rotation in \M{3} (elliptic)}
\label{rotation in M3 elliptic}
\end{figure}

The elliptic rotation is illustrated in Figure~\ref{rotation in M3 elliptic}, where 
 \(\mb{\Omega}=\tfrac{1}{\sqrt{8}}(-\e_{23}+3\e_{12})\), the axis of rotation, is proper.
The rotation of a point \(\tb{P}\) by the angle \(\theta=1.5\) is shown in Figure~\ref{rotation in M3 elliptic}(a)
and the rotation of a line \(\mb{\Lambda}\) by the same angle is shown in Figure~\ref{rotation in M3 elliptic}(b).
The rotation of \(\tb{P}\) is contained in the plane \(\mb{\Omega}\cdot\tb{P}\), which passes through  \(\tb{P}\)
and is perpendicular to \(\mb{\Omega}\). 
This plane is improper when \(\mb{\Omega}\) is proper.
Under the action of \(S\) with variable \(\theta\) ranging from \(0\) to \(2\pi\),
\(\tb{P}\) follows a trajectory in \M{2} shown in Figure~\ref{rotation in M3 elliptic} with the dotted curve.
It consists of points \(\tb{X}\) which satisfy
\(\norm{\tb{X}\vee\tb{C}}=\norm{\tb{P}\vee\tb{C}}\), where \(\tb{C}=\mb{\Omega}\wedge(\mb{\Omega}\cdot\tb{P})\)
 and \(\tb{X}\), \(\tb{P}\), and \(\tb{C}\) are assumed to be normalised.
The trajectory is an ellipse when \(\mb{\Omega}\) is proper, so the rotation is called elliptic.
Under elliptic rotation, the other points in the plane \(\mb{\Omega}\cdot\tb{P}\) follow elliptic trajectories centred on the point \(\tb{C}\) too.

A proper sphere of radius \(r\) centred at \(\tb{C}\) 
consists of normalised points \(\tb{X}\) for which \(\norm{\tb{X}\vee\tb{C}}=r\) and \(\tb{X}\vee\tb{C}\) is proper (\(\tb{C}\) is assumed to be normalised).
The points which constitute a proper sphere of radius \(r\) are at the distance \(r\) from the centre.
An improper sphere is defined in the same way, except that \(\tb{X}\vee\tb{C}\) is improper.
A proper sphere coincides with a hyperboloid of two sheets and an improper sphere coincides with a hyperboloid of one sheet;
the hyperboloids share the same axis which is parallel to the \(t\)-axis.
The elliptic rotation of a point \(\tb{P}\) around a  proper line \(\mb{\Omega}\) follows an ellipse,
where the plane  \(\mb{\Omega}\cdot\tb{P}\) intersect an improper sphere centered at 
\(\tb{C}=\mb{\Omega}\wedge(\mb{\Omega}\cdot\tb{P})\)
with the radius given by \(\norm{\tb{P}\vee\tb{C}}\).
The plane  \(\mb{\Omega}\cdot\tb{P}\) does not intersect the accompanying proper sphere as long as \(\mb{\Omega}\) is proper.

\begin{figure}[t!]
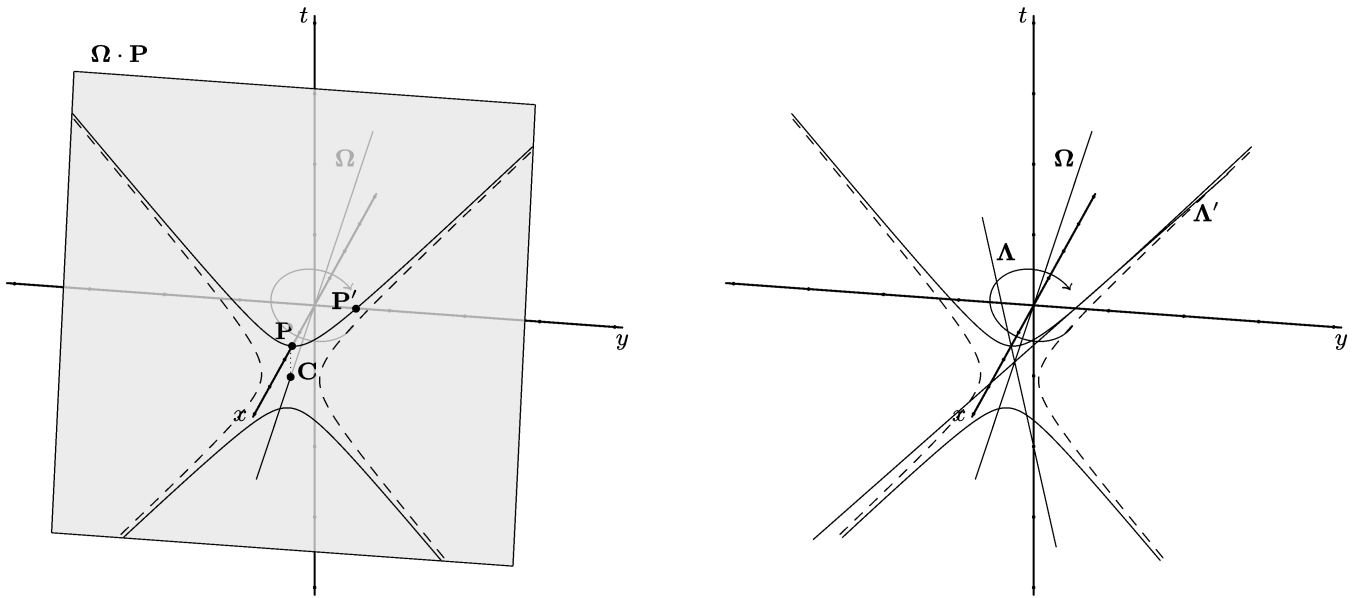
\hspace{-1cm}
\begin{subfloatenv}{}
\begin{asy}
import Figure3D;
metric = Metric(Minkowski);
Figure f = Figure(metric);

var P = Point(1,1.5,0,0);
f.point(P, "$\textbf{P}$",align=(0,-0.5,1),draw_orientation=false);

MV L = Line(0,0,0,-1,0,0.25);
L = L/norm(L);
f.plane(dot(L,P),"$\mathbf{\Omega}\cdot\textbf{P}$", align=(0,2.5,2),draw_orientation=false);
f.line(L, size=4,"$\mathbf{\Omega}$", position=0.9, align=(0,-1,0.5));

//int n = 10; for(int i: sequence(n)) { real theta = 2*i/n; write(theta); MV S = exp(-1/2*theta*L); f.point(S*P/S,string(i), draw_helper_lines=false, draw_orientation=false); }

real theta = 1.5;
MV S = exp(-1/2*theta*L); 
f.point(S*P/S,"$\textbf{P}'$", draw_helper_lines=false, draw_orientation=false,align=(1,-0.5,1));

//triple f(real theta) { MV S = exp(-1/2*theta*L); return totriple(S*P/S); }
triple f(real theta) { MV S = cosh(-1/2*theta)+L*sinh(-1/2*theta); return totriple(S*P/S); }
real x(real theta) { return f(theta).x; }
real y(real theta) { return f(theta).y; }
real z(real theta) { return f(theta).z; }
path3 p = graph(x, y, z, -2.765, 2.765, operator ..);
draw(p);

MV C = wedge(L,dot(L,P));
f.point(C, "$\textbf{C}$",align=(1,1.5,1),draw_orientation=false);
MV P1 = C*P/C;
//triple f(real theta) { MV S = exp(-1/2*theta*L); return totriple(S*P1/S); }
triple f(real theta) { MV S = cosh(-1/2*theta)+L*sinh(-1/2*theta); return totriple(S*P1/S); }
real x(real theta) { return f(theta).x; }
real y(real theta) { return f(theta).y; }
real z(real theta) { return f(theta).z; }
path3 p = graph(x, y, z, -2.4, 2.4, operator ..);
draw(p);

MV a = join(P,L);
a = a/norm(a);
MV T = exp(-1/2*norm(join(P,C))*wedge(e_0,a));
MV Q = T*C/T;
MV Q1 = C*Q/C;

//triple f(real theta) { MV S = exp(-1/2*theta*L); return totriple(S*Q/S); }
triple f(real theta) { MV S = cosh(-1/2*theta)+L*sinh(-1/2*theta); return totriple(S*Q/S); }
real x(real theta) { return f(theta).x; }
real y(real theta) { return f(theta).y; }
real z(real theta) { return f(theta).z; }
path3 p = graph(x, y, z, -2.4, 2.75, operator ..);
draw(p,dashed);

//triple f(real theta) { MV S = exp(-1/2*theta*L); return totriple(S*Q1/S); }
triple f(real theta) { MV S = cosh(-1/2*theta)+L*sinh(-1/2*theta); return totriple(S*Q1/S); }
real x(real theta) { return f(theta).x; }
real y(real theta) { return f(theta).y; }
real z(real theta) { return f(theta).z; }
path3 p = graph(x, y, z, -2.75, 2.4, operator ..);
draw(p,dashed);

\end{asy}
\end{subfloatenv}\hfill%
\begin{subfloatenv}{}
\begin{asy}
import Figure3D;
metric = Metric(Minkowski);
Figure f = Figure(metric);

var P = Point(1,1.5,0,0);

MV L = Line(0,0,0,-1,0,0.25);
L = L/norm(L);
f.line(L, size=4,"$\mathbf{\Omega}$", position=0.9, align=(0,-1,0.5));

var M = join(Point(1,1.5,0,0), Point(1,0,0,-2))/12;

f.line(M,"$\mathbf{\Lambda}$",position=0.9,draw_orientation=false);

real theta = 1.5;
MV S = exp(-1/2*theta*L); 
f.line(S*M/S,"$\mathbf{\Lambda}'$",position=0.9,draw_orientation=false);

MV P = Point(1,1.5,0,0);
//triple f(real theta) { MV S = exp(-1/2*theta*L); return totriple(S*P/S); }
triple f(real theta) { MV S = cosh(-1/2*theta)+L*sinh(-1/2*theta); return totriple(S*P/S); }
real x(real theta) { return f(theta).x; }
real y(real theta) { return f(theta).y; }
real z(real theta) { return f(theta).z; }
path3 p = graph(x, y, z, -2.765, 2.765, operator ..);
draw(p);

MV C = wedge(L,dot(L,P));
MV Q = C*P/C;
//triple f(real theta) { MV S = exp(-1/2*theta*L); return totriple(S*Q/S); }
triple f(real theta) { MV S = cosh(-1/2*theta)+L*sinh(-1/2*theta); return totriple(S*Q/S); }
real x(real theta) { return f(theta).x; }
real y(real theta) { return f(theta).y; }
real z(real theta) { return f(theta).z; }
path3 p = graph(x, y, z, -2.4, 2.4, operator ..);
draw(p);

MV a = join(P,L);
a = a/norm(a);
MV T = exp(-1/2*norm(join(P,C))*wedge(e_0,a));
MV Q = T*C/T;
MV Q1 = C*Q/C;

//triple f(real theta) { MV S = exp(-1/2*theta*L); return totriple(S*Q/S); }
triple f(real theta) { MV S = cosh(-1/2*theta)+L*sinh(-1/2*theta); return totriple(S*Q/S); }
real x(real theta) { return f(theta).x; }
real y(real theta) { return f(theta).y; }
real z(real theta) { return f(theta).z; }
path3 p = graph(x, y, z, -2.4, 2.75, operator ..);
draw(p,dashed);

//triple f(real theta) { MV S = exp(-1/2*theta*L); return totriple(S*Q1/S); }
triple f(real theta) { MV S = cosh(-1/2*theta)+L*sinh(-1/2*theta); return totriple(S*Q1/S); }
real x(real theta) { return f(theta).x; }
real y(real theta) { return f(theta).y; }
real z(real theta) { return f(theta).z; }
path3 p = graph(x, y, z, -2.75, 2.4, operator ..);
draw(p,dashed);

\end{asy}
\end{subfloatenv}
\caption{Rotation in \M{3} (hyperbolic)}
\label{rotation in M3 hyperbolic}
\end{figure}

The hyperbolic rotation of a point \(\tb{P}\) and a line \(\mb{\Lambda}\) is shown in Figure~\ref{rotation in M3 hyperbolic}, where
\(\mb{\Omega}=\tfrac{4}{\sqrt{15}}(-\e_{23}+\tfrac{1}{4}\e_{12})\) is improper and \(\theta=1.5\) as above.
Like the elliptic rotation, it is confined to the plane \(\mb{\Omega}\cdot\tb{P}\), but this plane is proper when \(\mb{\Omega}\) is improper.
The plane \(\mb{\Omega}\cdot\tb{P}\) intersects a proper sphere, 
which is centred at \(\tb{C}=\mb{\Omega}\wedge(\mb{\Omega}\cdot\tb{P})\) and whose radius equals  \(\norm{\tb{P}\vee\tb{C}}\),
along the hyperbola shown in Figure~\ref{rotation in M3 hyperbolic} with solid curves.
The same plane intersects an improper sphere of the same radius and also centred at \(\tb{C}\)
along the hyperbola shown with dashed curves.
As long as \(\mb{\Omega}\) is improper, the plane \(\mb{\Omega}\cdot\tb{P}\) will intersect
both the proper sphere and the improper one.
Under the action of \(S\) with variable \(\theta\) ranging from \(-\infty\) to \(+\infty\),
the point \(\tb{P}\) follows the upper branch of the hyperbola shown as a solid curve passing through \(\tb{P}\) in Figure~\ref{rotation in M3 hyperbolic}(a).
The effect of the action of \(S\) on the other points located in the plane \(\mb{\Omega}\cdot\tb{P}\) is reminiscent
of the action generated by a rotation in \M{2}, which is shown in Figure~\ref{rotation and translation M2}(a).

\begin{figure}[t]
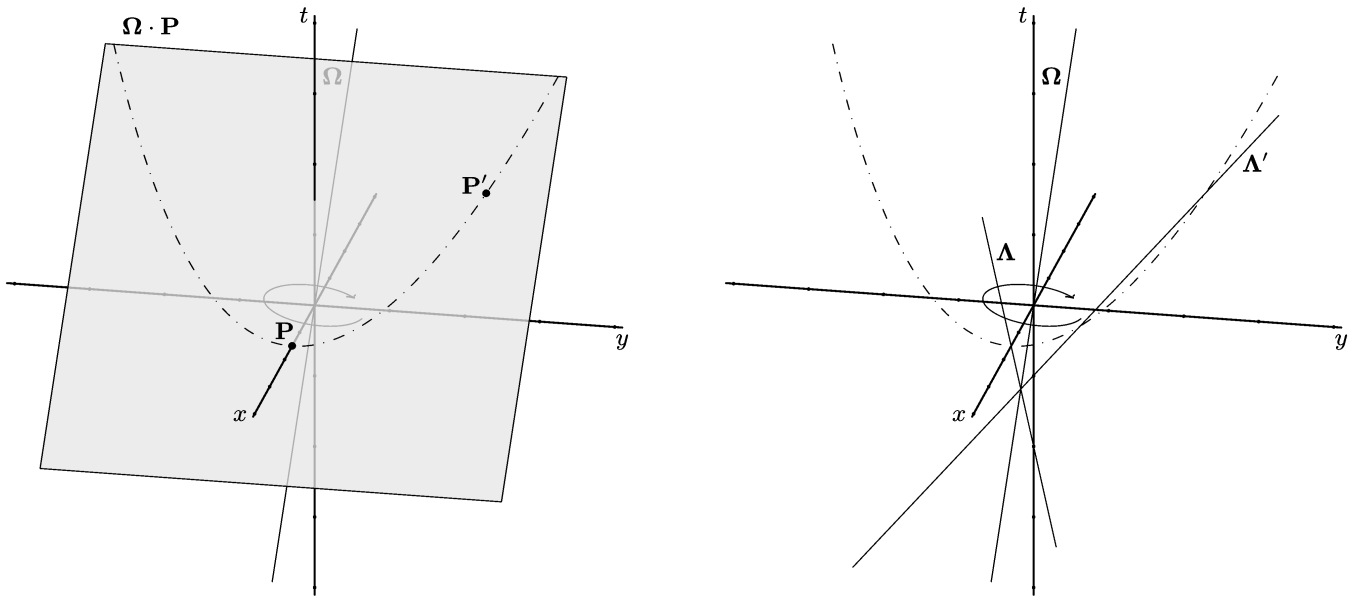
\hspace{-1cm}
\begin{subfloatenv}{}
\begin{asy}
import Figure3D;
metric = Metric(Minkowski);
Figure f = Figure(metric);

var P = Point(1,1.5,0,0);
f.point(P, "$\textbf{P}$",align=(0,-0.5,1),draw_orientation=false);

MV L = Line(0,0,0,-1,0,1);
f.plane(dot(L,P),"$\mathbf{\Omega}\cdot\textbf{P}$", align=(0,2.5,2),draw_orientation=false);
f.line(L, size=4,"$\mathbf{\Omega}$", position=0.9, align=(0,-1,0.5));

//int n = 10; for(int i: sequence(n)) { real theta = 2*i/n; write(theta); MV S = exp(-1/2*theta*L); f.point(S*P/S,string(i), draw_helper_lines=false, draw_orientation=false); }

real theta = 1.5;
MV S = exp(-1/2*theta*L); 
f.point(S*P/S,"$\textbf{P}'$", draw_helper_lines=false, draw_orientation=false,align=(1,-0.5,1));

//triple f(real theta) { MV S = exp(-1/2*theta*L); return totriple(S*P/S); }
triple f(real theta) { MV S = 1-1/2*theta*L; return totriple(S*P/S); }
real x(real theta) { return f(theta).x; }
real y(real theta) { return f(theta).y; }
real z(real theta) { return f(theta).z; }
path3 p = graph(x, y, z, -1.975, 1.975, operator ..);
draw(p,dashdotted);

\end{asy}
\end{subfloatenv}\hfill%
\begin{subfloatenv}{}
\begin{asy}
import Figure3D;
metric = Metric(Minkowski);
Figure f = Figure(metric);

var P = Point(1,1.5,0,0);
MV L = Line(0,0,0,-1,0,1);
f.line(L, size=4,"$\mathbf{\Omega}$", position=0.9, align=(0,-1,0.5));

var M = join(Point(1,1.5,0,0), Point(1,0,0,-2))/12;

f.line(M,"$\mathbf{\Lambda}$",position=0.9,draw_orientation=false);

real theta = 1.5;
MV S = exp(-1/2*theta*L); 
f.line(S*M/S,"$\mathbf{\Lambda}'$",position=0.9,draw_orientation=false);

//triple f(real theta) { MV S = exp(-1/2*theta*L); return totriple(S*P/S); }
triple f(real theta) { MV S = 1-1/2*theta*L; return totriple(S*P/S); }
real x(real theta) { return f(theta).x; }
real y(real theta) { return f(theta).y; }
real z(real theta) { return f(theta).z; }
path3 p = graph(x, y, z, -1.975, 1.975, operator ..);
draw(p,dashdotted);

\end{asy}
\end{subfloatenv}
\caption{Rotation in \M{3} (parabolic)}
\label{rotation in M3 parabolic}
\end{figure}

The parabolic rotation is a special case, since \(\tb{C}=\mb{\Omega}\wedge(\mb{\Omega}\cdot\tb{P})\)
is a point at infinity and \(\norm{\tb{X}\vee\tb{C}}=0\) for any point \(\tb{X}\) in \M{3}.
It is illustrated in Figure~\ref{rotation in M3 parabolic}, where \(\mb{\Omega}=-\e_{23}+\e_{12}\) is null,
\(\tb{P}=\e_{123}+\tfrac{3}{2}\e_{320}\), and \(\theta=1.5\).
More generally, the point \(\tb{P}\) shifts to \(\tb{P}(\theta)\) under the action of 
\(S=e^{-\tfrac{1}{2}\theta(-\e_{23}+\e_{12})}\), which depends on \(\theta\).
Since the bivector \(-\e_{23}+\e_{12}\) is null, 
\(S=1-\tfrac{1}{2}\theta(-\e_{23}+\e_{12})\) and 
\[\tb{P}(\theta) = S\tb{P}S^{-1}=
\big(1-\tfrac{1}{2}\theta(-\e_{23}+\e_{12})\big)
\big(\e_{123}+\tfrac{3}{2}\e_{320}\big)
\big(1+\tfrac{1}{2}\theta(-\e_{23}+\e_{12})\big),
\]
which gives
\[
\left\{
\begin{aligned}
&x(\theta)=\tfrac{3}{4}(2-\theta^2),\\
&y(\theta)=\tfrac{3}{2}\theta,\\
&t(\theta)=\tfrac{3}{4}\theta^2,
\end{aligned}
\right.
\]
for the coordinates of \(\tb{P}(\theta)\).
So, 
 \(\tb{P}(\theta)\) follows a parabolic trajectory (dash-dotted curve) and, therefore, the rotation is called parabolic.

Unlike Euclidean space, it is not possible to define a universal angular measure in \M{3}.
The parameter \(\theta\) in (\ref{boost and Euclidean rotation in M3}) can be used in a limited way to 
measure angles between planes passing through \(\mb{\Omega}\), the axis of rotation.
For instance, if \(\mb{\Omega}\) is parallel to the \(t\)-axis, then \(\theta\) is purely Euclidean as in \E{2}.
On the other hand, if \(\mb{\Omega}\) is parallel to the \(xy\)-plane, then \(\theta\) is purely pseudo-Euclidean as in \M{2}.

A composition of two boosts is a boost too only if the lines defining the boosts commute.
Let \(B_1=e^{-\tfrac{1}{2}\phi_1\mb{\Omega}_1}\) and \(B_2=e^{-\tfrac{1}{2}\phi_2\mb{\Omega}_2}\),
where \(\mb{\Omega}_1\) and \(\mb{\Omega}_2\) are finite lines in the \(xy\)-plane.
If the lines \(\mb{\Omega}_1\) and \(\mb{\Omega}_2\) commute, then \(S=B_1B_2\) is a boost.
Indeed, I can write \(B_1B_2=e^{-\tfrac{1}{2}\phi_1\mb{\Omega}_1}e^{-\tfrac{1}{2}\phi_2\mb{\Omega}_2}
=e^{-\tfrac{1}{2}(\phi_1\mb{\Omega}_1+\phi_2\mb{\Omega}_2)}\),
where \(\phi_1\mb{\Omega}_1+\phi_2\mb{\Omega}_2\) is a finite line in the \(xy\)-plane
and, therefore, \(S\) can be written in the exponential form for some line \(\mb{\Omega}\) in the \(xy\)-plane.
In other words, \(S\) is a boost as well.
Two finite lines in the \(xy\)-plane commute only if their attitudes coincide.
So, a composition of two boosts is a boost only if \(\mb{\Omega}_1=\pm\mb{\Omega}_2\) since
both lines are assumed to be normalised.
Otherwise, the composition will consist of a combination of a boost and either a Euclidean rotation or a translation in the \(xy\)-plane.

The action \(T\tb{P}T^{-1}\) of a spinor
\begin{equation}
T = e^{-\tfrac{1}{2}\lambda\e_0\wedge\tb{a}}=1-\tfrac{1}{2}\lambda\e_0\wedge\tb{a}
\end{equation}
on a point \(\tb{P}\) yields a translation of \(\tb{P}\) 
along a line perpendicular to the plane \(\tb{a}\) (if \(\tb{a}\) is null, then the translation is along a null line in \(\tb{a}\)).
The distance of the translation is defined if \(\tb{a}\) is improper and is given by \(\lambda\) if \(\tb{a}\) is normalised.
If \(\tb{P}'=T\tb{P}T^{-1}\) and \(\tb{a}=a\e_1+b\e_2+h\e_3\), then
\begin{equation}
\left\{
\begin{aligned}
&x'=x + \lambda a,\\
&y'=y+\lambda b\\
&t'=t-\lambda h.
\end{aligned}
\right.
\label{Translation 3d}
\end{equation}
So, a translation in \M{3} is a simple generalisation of a translation in \M{2}.

Any spinor in \M{3} can be written as  either
\begin{equation}
S_1=e^{-\tfrac{1}{2}(\theta-\lambda\I)\mb{\Lambda}},
\end{equation}
where \(\mb{\Lambda}\) is normalised (\(\mb{\Lambda}\) is simple and not null),
or
\begin{equation}
S_2=e^{-\tfrac{1}{2}\theta\mb{\Lambda}},
\end{equation}
where \(\mb{\Lambda}\) is null (\(\mb{\Lambda}\) can be simple or not simple).
The spinor \(S_1\) gives a proper motion consisting of a rotation around \(\mb{\Lambda}\) by the angle \(\theta\)
and a translation by \(\lambda\) along \(\mb{\Lambda}\).
The meaning of the angle \(\theta\) and the character of the rotation depend on whether \(\mb{\Lambda}\)
is proper or improper.
The action of \(S_2\) depends on whether it is simple or non-simple.
If \(\mb{\Lambda}\) is a null line (i.e.\ it is simple in addition to \(\mb{\Lambda}\cdot\mb{\Lambda}=0\)),
\(S_2\) gives  a proper motion consisting of a parabolic rotation in a null plane parallel to \(\mb{\Lambda}\).
Since \(\mb{\Lambda}\) is null in this case, it cannot be normalised and
the exact meaning of  \(\theta\) depends on the specific choice of the weight of \(\mb{\Lambda}\).
For instance, it might be convenient to change the weight such that the \(\e_{12}\) component
of \(\mb{\Lambda}\) equals unity, i.e.\ \(p_{12}=1\).
If \(\mb{\Lambda}\) is null and non-simple (i.e.\ \(\mb{\Lambda}\cdot\mb{\Lambda}=0\) and
\(\mb{\Lambda}\wedge\mb{\Lambda}\ne0\), which implies \(\mb{\Lambda}^3=\mb{\Lambda}\cdot(\mb{\Lambda}\wedge\mb{\Lambda})\ne0\)), 
the action of \(S_2\) is irreducible, i.e.\ it cannot be naturally and unambiguously reduced to more elementary actions.
Its action on a point \(\tb{P}\) is given by
\begin{equation}
S_2\tb{P}S_{2}^{-1}=
\tb{P} +\theta\tb{P}\times\mb{\Lambda}-\tfrac{1}{4}\theta^2\mb{\Lambda}\tb{P}\mb{\Lambda}
-\tfrac{1}{12}\theta^3 \tb{P}\times\mb{\Lambda}^3,
\end{equation}
which is not parabolic since it depends on \(\theta^3\).

In summary, proper motions in \M{3} consist of 
1) three kinds of rotation (elliptic, parabolic, and hyperbolic) including
boosts and Euclidean rotations,
2) translation (along proper, improper, or null directions),
3) a combination of elliptic or hyperbolic rotation and translation along the axis of rotation,
4) an irreducible action generated by a non-simple null bivector.

\section{Minkowski hyperspace \M{4} ( $3+1$ spacetime)}

For a hyperplane \(\tb{a}\), a plane \(\mb{\sigma}\), a line \(\mb{\Lambda}\), and a point \(\tb{P}\),
I have
\begin{equation}
\begin{split}
&\tb{a}^2=a^2+b^2+c^2-h^2,\quad\\
&\mb{\sigma}^2=p_{41}^2+p_{42}^2+p_{43}^2    -p_{23}^2-p_{31}^2-p_{12}^2  ,\quad\\
&\mb{\Lambda}^2=s_{234}^2+s_{314}^2+s_{124}^2-s_{321}^2,\quad
\end{split}
\end{equation}
and
\begin{equation}
\begin{split}
&\norm{\tb{a}}=|a^2+b^2+c^2-h^2|^{\tfrac{1}{2}},\quad
\norm{\mb{\sigma}}=|p_{23}^2+p_{31}^2+p_{12}^2-p_{41}^2-p_{42}^2-p_{43}^2|^{\tfrac{1}{2}},\quad\\
&\norm{\mb{\Lambda}}=|s_{234}^2+s_{314}^2+s_{124}^2-s_{321}^2|^{\tfrac{1}{2}},\quad
\norm{\tb{P}}=|w|.
\end{split}
\end{equation}
Unlike the 2- and 3-dimensional cases, % where \(\reverse{\tb{P}}=-\tb{P}\), 
 \(\tb{P}^2=-w^2\) in \M{4}.% since  \(\reverse{\tb{P}}=\tb{P}\).
% and the definition of the metric (note \(\tb{P}^2=w^2\) in \E{4}).

%A finite blade \(A_k\) is proper if \(A_k\reverse{A_k}>0\), improper if \(A_k\reverse{A_k}<0\),
%and null if \(A_k\reverse{A_k}=0\).
A hyperplane \(\tb{a}\) is proper if \(\tb{a}^2>0\),
a plane \(\mb{\sigma}\) is proper if \(\mb{\sigma}^2<0\),
and a line \(\mb{\Lambda}\) is proper if \(\mb{\Lambda}^2<0\).
Reversing the sign of the inequalities gives improper objects, while equality to zero for a finite object indicates that it is null.
Proper lines correspond to worldlines of objects moving at the speed less than the speed of light.
Null lines passing through a given point form a pair of hypercones.
Null planes and hyperplanes passing through the same point are tangent to the hypercones.

There are two ways of understanding geometric structures in \M{4}.
Firstly, the geometric intuitions developed in \M{2} and \M{3} can be applied in \M{4},
since many geometric structures in \M{4} can be reduced to those in \M{2} and \M{3}.
Secondly, many structures in \M{4} are essentially equivalent to those in \E{4}
with the suitable modification of the notion of perpendicularity.
For instance, \(\tb{a}\cdot\tb{P}\) is a line that passes through  \(\tb{P}\) and is perpendicular to  \(\tb{a}\).
This statement is essentially equivalent to the one I made in the section on \E{4} in \cite{sokolov2013clifford};
the difference is that the notion of perpendicularity in \M{4} is not the same as in \E{4}.
In general, the distance measure in \M{4} is an extension of that in \M{2} and \M{3}.
A limited angular measure can be defined for some hyperplanes, but it depends on the axis of rotation.
The other notable distinction between \M{4} and \E{4} is the presence of null geometric objects in \M{4},
which are not invertible and perpendicular to themselves.

For a non-simple trivector \(\mb{\Phi}\),
 the same decomposition that was used in \E{4} applies in \M{4}, 
provided that \(\mb{\Phi}\cdot\mb{\Phi}\ne0\). 
Namely, \(\mb{\Phi}=\mb{\Phi}_{fc}+\mb{\Phi}_{ic}\), where
\begin{equation}
\mb{\Phi}_{fc}=\left(1-\frac{(\mb{\Phi}\vee\mb{\Phi})\I}{2\mb{\Phi}\cdot\mb{\Phi}}\right)\mb{\Phi} \textrm{ and }
\mb{\Phi}_{ic}=\frac{(\mb{\Phi}\vee\mb{\Phi})\I}{2\mb{\Phi}\cdot\mb{\Phi}}\mb{\Phi}.
\end{equation}
In this decomposition, \(\mb{\Phi}_{fc}\) is a finite line and \(\mb{\Phi}_{ic}\) is a line at infinity.
%Planes that belong to the stack that represents \(\mb{\Phi}_{ic}\) are perpendicular to \(\mb{\Phi}_{fc}\).
%Furthermore, the components \(\mb{\Phi}_{fc}\) and \(\mb{\Phi}_{ic}\) commute.
If \(\mb{\Phi}\) is not simple and  \(\mb{\Phi}\cdot\mb{\Phi}=0\),  there are two possibilities.
Firstly, \(\mb{\Phi}=\e_0\wedge\mb{\pi}\) where the bivector \(\mb{\pi}\) is not simple.
The decomposition of \(\mb{\Phi}\) depends on whether \(\mb{\pi}\) can be decomposed into complementary planes.
For instance, \(\mb{\Phi}=\e_{310} +\e_{420}\) splits into two lines at infinity
\(\mb{\Phi}_1=\e_{310}\) and \(\mb{\Phi}_2=\e_{420}\).
Secondly,  \(\mb{\Phi}\) is null (i.e.\ its component \(\mb{\Phi}_0\), which passes through the origin, is null), 
in which case the decomposition is impossible.
For instance, \(\mb{\Phi}=\e_{321}+\e_{234}+\e_{310}\) cannot be decomposed (note that \(\mb{\Phi}^3\ne0\) in this case).

Let \(\mb{\pi}\) be a non-simple bivector, i.e.\ \(\mb{\pi}\wedge\mb{\pi}\ne0\).
If \(\mb{\pi}\) can be decomposed into two complementary planes \(\mb{\pi}_1\) and \(\mb{\pi}_2\),
then \(\mb{\pi}_1^2\) and \(\mb{\pi}_2^2\) satisfy
\begin{equation}
\mb{\pi}_{1,2}^4
-(\mb{\pi}\cdot\mb{\pi})\mb{\pi}_{1,2}^2
+\tfrac{1}{4}(\mb{\pi}\wedge\mb{\pi})^2=0,
\label{p1 p2 in M4}
\end{equation}
which is the same equation as in \E{4}.
But unlike \E{4}, \((\mb{\pi}\wedge\mb{\pi})^2\le0\) in \M{4} for any non-simple bivector
and, therefore, equation~(\ref{p1 p2 in M4}) has two solutions
\begin{equation}
\mb{\pi}_{1,2}^2=
\tfrac{1}{2}\left(
\mb{\pi}\cdot\mb{\pi}
\pm
\sqrt{({\mb{\pi}\cdot\mb{\pi}})^2 - (\mb{\pi}\wedge\mb{\pi})^2}
\right)
\end{equation}
unless \(\mb{\pi}\cdot\mb{\pi}=0\)  and \((\mb{\pi}\wedge\mb{\pi})^2=0\).

If \((\mb{\pi}\wedge\mb{\pi})^2>0\), the decomposition consists of two finite planes,
the axes of \(\mb{\pi}\), given by
\begin{equation}
\mb{\pi}_1
=
\mb{\pi}/ (1+\tfrac{1}{2}\mb{\pi}\wedge\mb{\pi} / \mb{\pi}_1^2),
\quad
%\textrm{and}\quad
\mb{\pi}_2
=
\mb{\pi}/(1+\tfrac{1}{2}\mb{\pi}\wedge\mb{\pi} / \mb{\pi}_2^2)
%\label{s1 and s2}
\end{equation}
or, alternatively,
\begin{equation}
\mb{\pi}_1
=
\frac{(1-\frac{1}{2}\mb{\pi}\wedge\mb{\pi} / \mb{\pi}_1^2)\mb{\pi}}
{1-\frac{1}{4}(\mb{\pi}\wedge\mb{\pi})^2 / \mb{\pi}_1^4},\quad
%\textrm{and}\quad
\mb{\pi}_2
=
\frac{(1-\frac{1}{2}\mb{\pi}\wedge\mb{\pi} / \mb{\pi}_2^2)\mb{\pi}}
{1-\frac{1}{4}(\mb{\pi}\wedge\mb{\pi})^2 / \mb{\pi}_2^4}.
\end{equation}
In this decomposition, neither \(\mb{\pi}_1\) nor \(\mb{\pi}_2\) is null.
One of the axes is proper and the other is improper.
There is no restriction on the value of \(\mb{\pi}\cdot\mb{\pi}\),
which can be positive, negative, or zero.
For instance, \(\mb{\pi}\cdot\mb{\pi}=0\) if \(\mb{\pi}=\e_{23}+\e_{41}\),
which still has a unique decomposition in \M{4}.

If \((\mb{\pi}\wedge\mb{\pi})^2=0\) and \(\mb{\pi}\cdot\mb{\pi}\ne0\),
the decomposition consists of a finite plane \(\mb{\pi}_1\) and a plane at infinity \(\mb{\pi}_2\) given by
\begin{equation}
\displaystyle\mb{\pi}_1=\left(1-\frac{\mb{\pi}\wedge\mb{\pi}}{2\mb{\pi}\cdot\mb{\pi}}\right)\mb{\pi},
\quad
\mb{\pi}_2=\frac{\mb{\pi}\wedge\mb{\pi}}{2\mb{\pi}\cdot\mb{\pi}}\mb{\pi}.
\end{equation}
In this decomposition, the finite axis \(\mb{\pi}_1\) cannot be null but it can be either proper or improper.

Let \(\mb{\pi}\) be a non-simple bivector
that satisfies \((\mb{\pi}\wedge\mb{\pi})^2=0\) and \(\mb{\pi}\cdot\mb{\pi}=0\).
In general,  \(\mb{\pi}\) cannot be decomposed into complementary axes,
e.g.\ \(\mb{\pi}=\e_{12}+\e_{41}+\e_{20}\) does not have a decomposition
(note that \(\mb{\pi}^3=\mb{\pi}\cdot(\mb{\pi}\wedge\mb{\pi})\ne0\) in this case).
If a decomposition exists, it may consist of two null planes or a null plane and a plane at infinity,
but the decomposition is not unique.
For example,  \(\mb{\pi}=-\e_{10}-\e_{23}+\e_{12}+\e_{41}+\e_{43}\)
can be decomposed into two null planes \(\mb{\pi}_1=(\e_4-\e_2)\wedge\e_3\)
and \(\mb{\pi}_2=(\e_4-\e_2+\e_0)\wedge\e_1\)
or into a plane at infinity \(\mb{\pi}'_1=\tfrac{1}{2}(-\e_1+\e_3)\wedge\e_{0}\)  and a null plane \(\mb{\pi}'_2=(\e_1+\e_3)\wedge(\e_2-\e_4-\tfrac{1}{2}\e_0)\).

Projection, rejection, reflection, and scaling are defined in \M{4} in the usual way.
The definition of the Spin group and its Lie algebra are the same as well.
Spinors are even and each spinor satisfies \(S\reverse{S}=1\).
Any proper motion in \M{4} can be generated by a spinor \(e^A\) where \(A\) is a bivector.

A normalised proper line passing through the origin can be parametrised by the angles \(\alpha\), \(\beta\), and \(\phi\) as follows:
\begin{equation}
\mb{\Lambda}=\e_{234}\cos{\alpha}\sin{\beta} \sinh{\phi}
+\e_{314}\sin{\alpha}\sin{\beta}\sinh{\phi}
+\e_{124}\cos{\beta}\sinh{\phi}
+\e_{321}\cosh{\phi}.
\end{equation}
It represents the worldline of an object moving in the direction \((\alpha,\beta)\) at the speed \(u=\tanh{\phi}\) if \(\phi>0\),
or in the opposite direction if \(\phi<0\).
The action of a spinor
\begin{equation}
B=e^{-\tfrac{1}{2}\phi_B\mb{\sigma}_B}, \textrm{ where}\quad
\mb{\sigma}_B=\e_{41}\cos{\alpha_B}\sin{\beta_B}+\e_{42}\sin{\alpha_B}\sin{\beta_B}+\e_{43}\cos{\beta_B}
\label{boost spinor M4}
\end{equation}
on the line \(\mb{\Lambda}\) yields a boost in the direction \((\alpha_B,\beta_B)\) if \(\phi_B>0\).
Note that the plane \(\mb{\sigma}_B\) is improper; it lies entirely in the hyperplane \(\e_4\), the \(xyz\)-hyperplane.
In the standard basis, the equation
\begin{equation}
\mb{\Lambda}'=B\mb{\Lambda}B^{-1},
\end{equation}
where \(\mb{\Lambda}'\) is parametrised by the angles  \(\alpha'\), \(\beta'\), and \(\phi'\),
corresponds to 
\begin{align*}
\cos{\alpha'}\sin{\beta'}\sinh{\phi'}
=\,
&\cos{\alpha}\sin{\beta}\sinh{\phi}
+\cos{\alpha_B}\sin{\beta_B}\cosh{\phi}\sinh{\phi_B}+\\*
&+\cos{\alpha}\sin{\beta}\sinh{\phi}\cos^2\!\alpha_B\sin^2\!\beta_B(\cosh{\phi_B}-1)+\\*
&+\sin{\alpha}\sin{\beta}\sinh{\phi}\sin{\alpha_B}\sin{\beta_B}\cos{\alpha_B}\sin{\beta_B}(\cosh{\phi_B}-1)+\\*
&+\cos{\beta}\sinh{\phi}\cos{\beta_B}\cos{\alpha_B}\sin{\beta_B}(\cosh{\phi_B}-1),\\
\sin{\alpha'}\sin{\beta'}\sinh{\phi'}
=\,
&\sin{\alpha}\sin{\beta}\sinh{\phi}
+\sin{\alpha_B}\sin{\beta_B}\cosh{\phi}\sinh{\phi_B}+\\*
&+\cos{\alpha}\sin{\beta}\sinh{\phi}\cos{\alpha_B}\sin{\beta_B}\sin{\alpha_B}\sin{\beta_B}(\cosh{\phi_B}-1)+\\*
&+\sin{\alpha}\sin{\beta}\sinh{\phi}\sin^2\!\alpha_B\sin^2\!\beta_B(\cosh{\phi_B}-1)+\\*
&+\cos{\beta}\sinh{\phi}\cos{\beta_B}\sin{\alpha_B}\sin{\beta_B}(\cosh{\phi_B}-1),\\
\cos{\beta'}\sinh{\phi'}
=\,
&\cos{\beta}\sinh{\phi}
+\cos{\beta_B}\cosh{\phi}\sinh{\phi_B}+\\*
&+\cos{\alpha}\sin{\beta}\sinh{\phi}\cos{\alpha_B}\sin{\beta_B}\cos{\beta_B}(\cosh{\phi_B}-1)+\\*
&+\sin{\alpha}\sin{\beta}\sinh{\phi}\sin{\alpha_B}\sin{\beta_B}\cos{\beta_B}(\cosh{\phi_B}-1)+\\*
&+\cos{\beta}\sinh{\phi}\cos^2\!\beta_B(\cosh{\phi_B}-1),\\
\cosh{\phi'}
=
&\cosh{\phi}\cosh{\phi_B}+\\*
+(\cos{\alpha}&\sin{\beta}\cos{\alpha_B}\sin{\beta_B}+
\sin{\alpha}\sin{\beta}\sin{\alpha_B}\sin{\beta_B}+\cos{\beta}\cos{\beta_B})\sinh{\phi}\sinh{\phi_B}.
\end{align*}
Diving the first three equations by the last and substituting
\begin{equation*}
\begin{split}
&(u_x,u_y,u_z)=(\cos{\alpha}\sin{\beta}\tanh{\phi},\sin{\alpha}\sin{\beta}\tanh{\phi},\cos{\beta}\tanh{\phi}),\\
&(u'_x,u'_y,u'_z)=(\cos{\alpha'}\sin{\beta'}\tanh{\phi'},\sin{\alpha'}\sin{\beta'}\tanh{\phi'},\cos{\beta'}\tanh{\phi'}),\\
&(v_x,v_y,v_z)=(\cos{\alpha_B}\sin{\beta_B}\tanh{\phi_B},\sin{\alpha_B}\sin{\beta_B}\tanh{\phi_B},\cos{\beta_B}\tanh{\phi_B})
\end{split}
\end{equation*}
yields
\begin{equation}
\left\{
\begin{aligned}
&u'_x
=
\frac{
{\gamma}^{-1}u_x
+v_x
+v_x(u_xv_x+u_yv_y+u_zv_z)\gamma(\gamma+1)^{-1}
}
{1+(u_xv_x+u_yv_y+u_zv_z)},\\
&u'_y
=
\frac{
{\gamma}^{-1}u_y
+v_y
+v_y(u_xv_x+u_yv_y+u_zv_z)\gamma(\gamma+1)^{-1}
}
{1+(u_xv_x+u_yv_y+u_zv_z)},\\
&u'_z
=
\frac{
{\gamma}^{-1}u_z
+v_z
+v_z(u_xv_x+u_yv_y+u_zv_z)\gamma(\gamma+1)^{-1}
}
{1+(u_xv_x+u_yv_y+u_zv_z)},\\
\end{aligned}
\right.
\label{velocity transform in M4 components}
\end{equation}
where  \(\gamma=\cosh{\phi_B}\) is the Lorentz factor corresponding to \(v=\sqrt{v_x^2+v_y^2+v_z^2}\).

The Lorentz transformation
\begin{equation}
\tb{P}'=B\tb{P}B^{-1}
\end{equation}
where \(\tb{P}=\e_{1234}+x\e_{2340}+y\e_{3140}+z\e_{1240}+t\e_{3210}\),
\(\tb{P}'=\e_{1234}+x'\e_{2340}+y'\e_{3140}+z'\e_{1240}+t'\e_{3210}\),
and the spinor \(B\) is given by (\ref{boost spinor M4}),
corresponds to 
\begin{equation}
\begin{aligned}
x'
=
x
&+t\cos{\alpha_B}\sin{\beta_B}\sinh{\phi_B}
+x\cos^2\!{\alpha_B}\sin^2\!{\beta_B}(\cosh{\phi_B}-1)+\\
&+y\cos{\alpha_B}\sin{\beta_B}\sin{\alpha_B}\sin{\beta_B}(\cosh{\phi_B}-1)
+z\cos{\alpha_B}\sin{\beta_B}\cos{\beta_B}(\cosh{\phi_B}-1),\\
y'
=
y
&+t\sin{\alpha_B}\sin{\beta_B}\sinh{\phi_B}
+x\sin{\alpha_B}\sin{\beta_B}\cos{\alpha_B}\sin{\beta_B}(\cosh{\phi_B}-1)+\\
&+y\sin^2\!{\alpha_B}\sin^2\!{\beta_B}(\cosh{\phi_B}-1)
+z\sin{\alpha_B}\sin{\beta_B}\cos{\beta_B}(\cosh{\phi_B}-1),\\
z'
=
z
&+t\cos{\beta_B}\sinh{\phi_B}
+x\cos{\beta_B}\cos{\alpha_B}\sin{\beta_B}(\cosh{\phi_B}-1)+\\
&+y\cos{\beta_B}\sin{\alpha_B}\sin{\beta_B}(\cosh{\phi_B}-1)
+z\cos^2\!{\beta_B}(\cosh{\phi_B}-1),\\
t'=t&\cosh{\phi_B}
+x\cos{\alpha_B}\sin{\beta_B}\sinh{\phi_B}
+y\sin{\alpha_B}\sin{\beta_B}\sinh{\phi_B}
+z\cos{\beta_B}\sinh{\phi_B}
\end{aligned}
\end{equation}
and
\begin{equation} 
\left\{
\begin{aligned}
&x'=x
+\gamma tv_x
+(\gamma-1)\frac{(xv_x+yv_y+zv_z)v_x}{v_x^2+v_y^2+v_z^2},\\
&y'=y
+\gamma tv_y
+(\gamma-1)\frac{(xv_x+yv_y+zv_z)v_y}{v_x^2+v_y^2+v_z^2}, \\
&z'=z
+\gamma tv_z
+(\gamma-1)\frac{(xv_x+yv_y+zv_z)v_z}{v_x^2+v_y^2+v_z^2}, \\
&t'=\gamma(t + xv_x+yv_y+zv_z).
\end{aligned}
\right.
\end{equation}

Consider a spinor
\begin{equation}
R=e^{-\tfrac{1}{2}\alpha\mb{\sigma}_E},\quad \textrm{where } 
\mb{\sigma}_E=\e_{23}\cos{\alpha_E}\sin{\beta_E}+\e_{31}\sin{\alpha_E}\sin{\beta_E}+\e_{12}\cos{\beta_E}.
\label{Euclidean spinor M4}
\end{equation}
The bivector  \(\mb{\sigma}_E\) dually represents a proper plane in \M{4}, which passes through the \(t\)-axis.
The action of \(R\) yields a rotation around this plane.
In the \(xyz\)-hyperplane, 
this rotation is equivalent to a 3-dimensional Euclidean rotation by the angle \(\alpha\) 
around the line \(\e_4\wedge\mb{\sigma}_E\), which lies in the \(xyz\)-hyperplane and points in the direction \((\alpha_E,\beta_E)\) within this hyperplane.
The same Euclidean rotation is induced by \(R\) in any other hyperplane parallel to \(\e_4\).
For an arbitrary finite line \(\mb{\Lambda}_E\), which lies in the \(xyz\)-hyperplane, 
a 3-dimensional Euclidean rotation around \(\mb{\Lambda}_E\)  can be generated by 
the spinor \(R\) with \(\mb{\sigma}_E=\e_4\cdot\mb{\Lambda}_E\).

Proper motions induced by the action \(SMS^{-1}\) of a spinor \(S\) on multivector \(M\)
can be summarised as follows:

\(S = e^{-\tfrac{1}{2}\lambda\e_0\wedge\tb{a}}\)
yields a translation along the lines perpendicular to \(\tb{a}\) 
(if \(\tb{a}\) is null, then the translation is along the null lines in \(\tb{a}\)).
The distance of translation is given by \(\lambda\) if \(\tb{a}\) is normalised and improper.

\(S=e^{-\tfrac{1}{2}\theta\mb{\sigma}}\),
where \(\mb{\sigma}\) is an arbitrary normalised or null plane,
yields a rotation around  \(\mb{\sigma}\) by the angle \(\theta\), which includes boosts and Euclidean rotations.
The character of rotation is elliptic, parabolic, or hyperbolic depending on whether the axis of rotation \(\mb{\sigma}\)
is proper, null, or improper.
The parameter \(\theta\) can be used as an angular measure for hyperplanes passing through \(\mb{\sigma}\).

\(S=e^{-\tfrac{1}{2}(\theta-\lambda\tb{a}\I)\mb{\sigma}}\),
where \(\mb{\sigma}\) is a normalised plane, \(\tb{a}\) is a normalised hyperplane, and \(\tb{a}\cdot\mb{\sigma}=0\),
yields a combination of an elliptic (or hyperbolic) rotation and a translation in the axis of rotation, along the lines perpendicular to \(\tb{a}\)
in the direction that depends on the orientation of \(\tb{a}\) and \(\mb{\sigma}\).

\(S=e^{-\tfrac{1}{2}(\theta_1\mb{\sigma}_1+\theta_2\mb{\sigma}_2)}\),
where \(\mb{\sigma}_1\) and \(\mb{\sigma}_2\) are normalised and complementary
(for definiteness, assume \(\mb{\sigma}_1\) is proper and \(\mb{\sigma}_2\) is improper),
yields a double rotation, i.e.\ an elliptic rotation by the angle \(\theta_1\) around \(\mb{\sigma}_1\)
and a hyperbolic rotation by angle \(\theta_2\) around \(\mb{\sigma}_2\).
As a special case, this includes a combination of a boost and a Euclidean rotation complementary to the boost,
which corresponds to a boost in the direction along a line and a 3-dimensional Euclidean rotation around the same line.
This type of double rotation is usually called loxodromic, as it produces loxodromic curves
around the axis of the Euclidean rotation.

\(S=e^{-\tfrac{1}{2}(\theta_1\mb{\eta}_1+\theta_2\mb{\eta}_2)}\),
where \(\mb{\eta}_1\) and \(\mb{\eta}_2\) are complementary planes
that are either both null or one null and the other is at infinity,
yields a transformation characterised by many invariant axes, either null or at infinity.
In this sense, it is reminiscent of the isoclinic rotation in \E{4}.

\(S=e^{-\tfrac{1}{2}\theta\mb{\eta}}\),
where \(\mb{\eta}\) is an irreducible non-simple null bivector 
(\(\mb{\eta}\cdot\mb{\eta}=0\), \((\mb{\eta}\wedge\mb{\eta})^2=0\), and \(\mb{\eta}^3\ne0\)),
yields a transformation that fails to reduce to either a double rotation in complementary null planes
or a combination of a rotation in a null plane and a translation in a complementary direction.
The action of \(S\) on a point \(\tb{P}\) can  be expressed as
\begin{equation}
S\tb{P}S^{-1}=
\tb{P} +\theta\tb{P}\times\mb{\eta}
-\tfrac{1}{4}\theta^2\mb{\eta}\tb{P}\mb{\eta}
-\tfrac{1}{12}\theta^3\tb{P}\times\mb{\eta}^3,
\end{equation}
so it does not generate a parabolic trajectory.

\bibliographystyle{plain}
\bibliography{g.bib}

\begin{thebibliography}{1}

\bibitem{gunn2011geometry}
Charles Gunn.
\newblock {\em Geometry, Kinematics, and Rigid Body Mechanics in Cayley-Klein
  Geometries}.
\newblock PhD thesis, TU Berlin, 2011.

\bibitem{sokolov2013clifford}
Andrey Sokolov.
\newblock Clifford algebra and the projective model of homogeneous metric
  spaces: Foundations, 2013.
\newblock arXiv:1307.2917 [math.MG].

\end{thebibliography}

\end{document}